\documentclass[final, 4p, 12pt]{elsarticle}
\usepackage{amssymb}
  \usepackage{amsmath}
  \usepackage{amsthm}
   \theoremstyle{plain}
  \usepackage{bm}
  
\newcommand\+{\mkern-4mu}
  \usepackage{verbatim}
  \usepackage{textcomp} 
\usepackage{booktabs} 
\usepackage{dcolumn}  
\usepackage{multirow} 
\newcommand{\qqed}{\tag*{$\square$}}
\usepackage{lettrine}

\usepackage{algorithm}
\usepackage{algorithmicx}
\usepackage{algpseudocode}
\makeatletter \renewcommand*{\ALG@name}{Method}
\usepackage{mathrsfs}  

\usepackage{simplemargins}
\setallmargins{.8in}
\newtheorem{thm}{Theorem}
\makeatletter
\renewcommand*\env@matrix[1][c]{\hskip -\arraycolsep
  \let\@ifnextchar\new@ifnextchar
  \array{*\c@MaxMatrixCols #1}}
\makeatother

\usepackage{url} 
\newtheorem{cor}{Corollary}
\newtheorem{ex}{Example}
\newtheorem{lem}{Lemma}

\newcommand{\ndiv}{\hspace{-4pt}\not|\hspace{2pt}}
\theoremstyle{definition}
\newtheorem{defn}{Definition}
\newtheorem*{pf}{Proof}

\DeclareFontFamily{U}{matha}{\hyphenchar\font45}
\DeclareFontShape{U}{matha}{m}{n}{
      <5> <6> <7> <8> <9> <10> gen * matha
      <10.95> matha10 <12> <14.4> <17.28> <20.74> <24.88> matha12
      }{}
\DeclareSymbolFont{matha}{U}{matha}{m}{n}
\DeclareFontFamily{U}{mathx}{\hyphenchar\font45}
\DeclareFontShape{U}{mathx}{m}{n}{
      <5> <6> <7> <8> <9> <10>
      <10.95> <12> <14.4> <17.28> <20.74> <24.88>
      mathx10
      }{}
\DeclareSymbolFont{mathx}{U}{mathx}{m}{n}

\DeclareMathSymbol{\obot}         {2}{matha}{"6B}
\DeclareMathSymbol{\bigobot}       {1}{mathx}{"CB}
\newcommand{\oant}{\mbox{OA}(N,k,2,s)}

\newcommand{\oand}{{\rm OA}(\lambda n^s,k,n,s)}
\newcommand{\oano}{{\rm OA}(n^s,k,n,s)}

\DeclareMathOperator{\Tr}{Tr}

\newcommand{\oal}{\mbox{OA}(\lambda 2^s,k,2,s)}

\newcommand{\I}{{\bf I}}
\newcommand{\J}{{\bf J}}
\newcommand{\sss}{\mathbf{s}}
\newcommand{\1}{{\bf 1}}
\newcommand{\0}{{\bf 0}}

\newcommand{\bbb}{{\bf{b}}}
\newcommand{\bb}{{\bf{b}}}

\newcommand{\p}{{\bf p}}

\newcommand{\x}{{\bf x}}

\newcommand{\vv}{{\bf v}}

\newcommand{\y}{{\bf y}}
\newcommand{\z}{{\bf z}}

\newcommand{\h}{{\bf h}}

\newcommand{\Y}{{\bf Y}}

\newcommand{\ZZ}{{\bf Z}}

\newcommand{\CCC}{\mathbb{C}}

\newcommand{\E}{{\bf E}}
\newcommand{\X}{{\bf X}}
\newcommand{\w}{{\bf w}}

\newcommand{\B}{{\bf B}}
\newcommand{\Gld}{G^{\rm LD}}

\newcommand{\ddd}{{\bf d}}
\newcommand{\D}{{\bf D}}
\newcommand{\CC}{{\bf C}}

\newcommand{\A}{{\bf{A}}}
\newcommand{\HH}{{\bf{H}}}

\newcommand{\etal}{\emph{et al.}}

\newcommand{\uuu}{{\bf u}}

\usepackage{hyperref}

\newcommand{\F}{\mathbb{F}}
\newcommand{\Q}{\mathbb{Q}}
\newcommand{\R}{\mathbb{R}}
\newcommand{\C}{\mathbb{C}}
\begin{document}
\begin{frontmatter}
\title{Finding the dimension of a non-empty  orthogonal array polytope}
\author[AFIT]{Dursun A.~Bulutoglu}
\ead{dursun.bulutoglu@gmail.com}
\address[AFIT]{Department of Mathematics and Statistics, Air Force Institute of Technology,\\Wright-Patterson Air Force Base, Ohio 45433, USA}
\journal{Discrete Optimization}
\begin{abstract}
By using representation theory, we reduce the size of the  set of possible values for the dimension of the convex hull of all feasible points of an orthogonal array (OA) defining  integer linear description (ILD). Our results 
address the conjecture that if this polytope is non-empty, then it is full-dimensional 
within the affine space where all the feasible points of the ILD's linear description (LD) relaxation   lie, raised by Appa \etal,\, [On multi-index assignment polytopes, Linear Algebra and its Applications 416 (2-3) (2006), 224--241]. In particular, our theoretical results
provide a sufficient condition for this polytope 
to be full-dimensional within the LD relaxation affine space 
when it is non-empty.
 This sufficient condition implies all the known non-trivial values of the dimension of the $(k,s)$ assignment polytope.
 However, our results suggest that the conjecture mentioned  above may not be true. More generally, we provide previously unknown restrictions on the feasible values of the dimension of the convex hull of all feasible points  of our OA defining ILD. We also determine all possible corresponding sets of equality constraints up to equivalence that can potentially be implied by the integrality constraints of this ILD. Moreover, we find additional restrictions on the dimension of the convex hull of all  feasible points, and larger sets of corresponding equality constraints for the $n=2$ and even $s$ cases. Each of these cases possesses symmetries that do not necessarily exist in the $3\leq n$ or odd $s$ cases.
 Finally, we discuss how to decrease  the number of possible values for the dimension of the convex hull of all feasible points of an arbitrary ILD as well as generate sets of corresponding equality constraints with the zero right hand side. These are the only sets of zero right hand side equality constraints up to equivalence that can potentially be implied by the integrality constraints of the ILD.
\end{abstract}
\end{frontmatter}
\section{Introduction}
An integer linear description (ILD) is a system of constraints of the form
 \begin{align}
& \quad  \A \x=\bb,\quad \B\x\leq  \ddd, \quad \x \in 
\mathbb{Z}^n, \label{eqn:geneqILP}
\end{align}
where $\A$ and $\B$ are  $m_1 \times n$ and $m_2 \times n$
constraint matrices, $\bb\in \mathbb{R}^{m_1}$,
 $\ddd \in \mathbb{R}^{m_2}$. 
 Let $P^{{\rm ILD}(\ref{eqn:geneqILP})}_{I}$ be the convex hull of all  feasible points of ILD~(\ref{eqn:geneqILP}). If $P^{{\rm ILD}(\ref{eqn:geneqILP})}_{I}$ is bounded  or $\bb \in \mathbb{Q}^{m_1}$, 
  $\ddd\in \mathbb{Q}^{m_2}$, and the matrices $\A$, $\B$ have only rational values, then $P^{{\rm ILD}(\ref{eqn:geneqILP})}_{I}$ is a polyhedron, and its facets are 
  its dim$(P^{{\rm ILD}(\ref{eqn:geneqILP})}_{I})-1$-dimensional faces.  Throughout the paper we assume that either $P^{{\rm ILD}(\ref{eqn:geneqILP})}_{I}$ is bounded  or $\bb \in \mathbb{Q}^{m_1}$, 
  $\ddd\in \mathbb{Q}^{m_2}$, and the matrices $\A$, $\B$ have only rational values.
 It is well-known that knowing facets of $P^{{\rm ILD}(\ref{eqn:geneqILP})}_{I}$  greatly decreases the time it takes to find a solution of ILD~(\ref{eqn:geneqILP})  or prove that no solution exists. However, determining whether a face of $P^{{\rm ILD}(\ref{eqn:geneqILP})}_{I}$ is a facet  
 requires knowing dim$(P^{{\rm ILD}(\ref{eqn:geneqILP})}_{I})$, and determining dim$(P^{{\rm ILD}(\ref{eqn:geneqILP})}_{I})$ is a difficult problem in its own right. 
 
Next we define  orthogonal arrays (OAs).
\begin{defn}\label{def:OA}
Let $\lambda \geq 1$, $n \geq 2$, $k \geq 1$ be integers, and $s$ be an integer such that  $1
\leq s \leq k$.  A $\lambda n^s \times k$ array $\Y$ whose entries are symbols from $\{l_1,\ldots,l_n\}$ is an {\em OA} of strength $s$, denoted by $\oand$, if each of the $n^s$ symbol
combinations from $\{l_1,\ldots,l_n\}^s$ appears $\lambda$ times in every $\lambda n^s \times  s $ subarray of  $\Y$.
 \end{defn}
 An OA($n^2,3,n,2$) is equivalent to an $n \times n$  {\em Latin square} and
 an OA($n^2,k,n,2$) is equivalent to $k-2$ 
 {\em mutually orthogonal $n \times n$ Latin squares}~\cite{Hedayat1999}.
 For $\lambda=1$, an $\oano$ is a {\em $(k,s)$ assignment 
 of order $n$}~\cite{Appa2006}.
  
Let  $\x \in \mathbb{Z}^{n^k}$, and $x(i_1,\ldots,i_k)$ 
  be the number of times the  symbol combination $(i_1,\ldots,i_k)$ such that  $(i_1,\ldots,i_k)^{\top} \in \{l_1,\ldots,l_n\}^k$ appears in an $\oand$.  Then, $\x$ is called the {\em frequency vector} of an $\oand$ and must be a feasible point of the ILD 
\begin{equation}\label{ilp:BF}
\begin{array}{ll}
  & \sum_{\{i_{1},\ldots,i_{k}\} \backslash \{i_{j_1}, \ldots,i_{j_s}\} \in \{l_1,\ldots,l_n\}^{k-s} } {x(i_1,\ldots,i_k)}= \lambda,\\
 &  0 \leq x(i_1,\ldots,i_k) \leq  p_{\max}, \quad x(i_1,\ldots,i_k) \in 
 \mathbb{Z}, \quad \text{for} \  (i_1,\ldots,i_k)^{\top} \in 
 \{l_1,\ldots,l_n\}^k, 
 \end{array}
\end{equation}
for each $\{j_1, \ldots,j_s\}\subseteq \{1,\ldots,k\}$ and each vector
 $(i_{j_1}, \ldots, i_{j_s})^{\top} \in \{l_1,\ldots,l_n\}^{s}$~\cite{Bulutoglu2008}, where $p_{\max}\leq \lambda$ is a positive integer computed as in~\cite{Bulutoglu2008}. For $\lambda=1$, ILD~(\ref{ilp:BF}) is the ILD formulation for the   {\em $(k,s)$ assignment problem of order $n$ ($(k,s)AP_n$)} in Appa~\etal~\cite{Appa2006,AppaJanssen2006,Appa2004}. 
For general $\lambda$, we call the constraint satisfaction problem  formulated by ILD~(\ref{ilp:BF}), the {\em $\oand$ problem}.

For $\lambda=1$, the convex hull of all the integer points satisfying ILD~(\ref{ilp:BF})
 is called the {\em $(k, s)$ assignment polytope}, denoted by {\em $P^{(k,s)}_{ n;I}$}~\cite{Appa2004}, and all the feasible points in $\mathbb{R}^{n^k}$ of the linear description (LD) relaxation of  ILD~(\ref{ilp:BF}) is
called the {\em linear $(k,s)$ assignment polytope},  denoted by 
$P_n^{(k,s)}$~\cite{Appa2004}.  For general $\lambda$, we call the corresponding concepts, {\em $(k, s, \lambda)$ orthogonal array polytope} denoted by $P^{(k,s,\lambda)}_{ n;I}$ and {\em $(k, s, \lambda)$ linear orthogonal array polytope} denoted by $P_n^{(k,s,\lambda)}$.  
In studying the facets  of $P_{n;I}^{(k,s)}$, 
Appa~\etal~{\cite{Appa2006}} tabulated Table~\ref{tab:PnI}  and conjectured that dim$(P_{n;I}^{(k,s)})=\text{dim}(P_{n}^{(k,s)})$
 provided that $P_{n;I}^{(k,s)}\neq \emptyset$.
In this paper, we address this conjecture by using representation theory. 

Throughout the paper, let $\delta_{(x,y)}=1$ if $x=y$, $\delta_{(x,y)}=0$ if $x\neq y$, and $E[a,b]$ be the set of even integers in the closed interval $[a,b]$. Then, the following  theorem is the main result of this paper.
\begin{thm}\label{thm:dimfin}
For the parameters $k,s,\lambda, n$ of  $P^{(k,s,\lambda)}_{n;I}$, let 
\begin{align*}
\Omega_1=\left\{d \in E[s+2,k]\mid\lambda {d-1\choose d-s-1} \equiv 0 \quad ({\rm mod}\  2)\right\}, \quad 	\Omega'_1=\begin{cases}\{k\} &\text{if $k$ is odd,}\\
			\emptyset &\text{otherwise,}
		\end{cases}
\end{align*}
$$\Omega_2=\left\{d \in  \{s+1,\ldots,k\} \mid \lambda \binom{d-1}{d-1-s}\equiv 0 \quad {\rm (mod } \ n)\right\}.$$ 
Then  $\exists  T_1\subseteq \Omega_1$, $\,T_1'\subseteq \Omega_1'$, and $T_2 \subseteq \Omega_2$ such that
	\begin{align*} 
		{\rm dim}(P^{(k,s,\lambda)}_{ n;I})=\begin{cases}
			2^k-\sum_{j=0}^s{k \choose j}-\left(\sum_{d\in T_1 }\left({k \choose d-1}+
			{k \choose d}\right)\right)-\delta_{(|T_1'|,1)} & \text{if $n=2$ and $s$ is even,} \\
			n^k-\sum_{j=0}^s{k \choose j}(n-1)^j-\sum_{d\in T_2 }{k \choose d}(n-1)^{d} & \text{otherwise.}
		\end{cases}
	\end{align*}	
\end{thm}
We will prove Theorem~\ref{thm:dimfin} by showing that the known symmetries of the feasible set of ILD~(\ref{ilp:BF}) drastically decrease the number of feasible values of dim$(P_{n;I}^{(k,s)})$, where a {\it symmetry of the feasible set} of an ILD is a permutation of its variables that sends a feasible point to a feasible point. The set of all symmetries 
of the feasible set of an ILD is called the {\em symmetry group} of the ILD.
Throughout the rest of the paper, we will be developing the needed theory to prove Theorem~\ref{thm:dimfin}.
\begin{table}[t]\label{tab:PnI}
\caption{Known values of dim($P^{(k,s)}_{n;I}$)} 
\begin{tabular}{l@{\hskip .4in} l@{\hskip .4in} l@{\hskip .4in} l}
\hline
$(k,s)$ & $n$ & dim($P^{(k,s)}_{n;I}$) & Reference\\
 \hline
 $(k, 0),\, \forall k \in \mathbb{Z}^{\geq 1} $   & $\geq 0$   &     $n^k-1$             & Appa~\etal~{\cite{Appa2006}}\\   
 $(2,1)$     & $\geq 2$  & $(n-1)^2$                 &  Balinski and Russakoff~\cite{Balinski1974}\\ 
 $(3,1)$     &  $\geq 3$ &   $n^3-3n+2$               &  
 Euler~\cite{Euler1987}, Balas and Saltzman~\cite{Balas1989}  \\ 
  $(3,2)$    & $\geq 4$  &    $(n-1)^3$              &   
  Euler~\etal~\cite{Euler1986}\\ 
   $(4,2)$    & $\geq 4,\, \neq 6$   &    $n^4 - 6n^2 + 8n -3$ &    
   Appa~\etal~\cite{AppaJanssen2006}\\ 
  $(k,k),\, \forall k \in \mathbb{Z}^{\geq 1}$    &  $\geq 0$ &    $0$              &  Appa~\etal~{\cite{Appa2006}}\\  
\hline
 \end{tabular}
\label{table1}
\end{table}

  A group $G$ with identity $e$ is said to {\em act} on a set $X$ if for each $(g,x) \in G\times X$, $gx \in X$, $ex=x$, and for each $g,h \in G$ we have $g(hx)=(gh)x$. Such a group action is called {\em transitive} if for each pair $(x_1,x_2)\in X\times X$, there exists $g \in G \, \ni \, gx_1=x_2$.
The following two definitions locate a subgroup
of the symmetry group of ILD~(\ref{ilp:BF}), and describe the action of this subgroup on the feasible points of ILD~(\ref{ilp:BF}).
 \begin{defn}
Two   $\oand$ are   \textit{isomorphic} if one can be obtained from the other by applying a sequence of permutations (including the identity) to the rows, columns,  and the elements of $\{l_1,\ldots,l_n\}$ 
within each column~\cite{Stufken2007}.
\end{defn}
\begin{defn}Let $\X$ be an $N$ row, $k$ column array with symbols from  $\{l_1,\ldots,l_n\}$. Then each of the $(n!)^kk!$ operations that involve permuting columns and the elements of $\{l_1,\ldots,l_n\}$ within each column of $\X$ is called an {\em isomorphism operation}. The set of all isomorphism operations forms a group $G^{\rm iso}(k,n)$ called the {\em paratopism group}~\cite{Egan2016}. 
\end{defn}
The group $G^{\rm iso}(k,n)$ acts on $\oand$, and
$G^{\rm iso}(k,n)$ is isomorphic to $S_n \wr S_k$~\cite{Egan2016}, where $S_n\wr S_k$ is the wreath product of the
 symmetric group  $S_n$ of degree $n$  and the symmetric group $S_k$ of degree $k$. The definition of the wreath product of groups can be found in~\cite{Rotman1994}.
 
 The {\it symmetry group $\Gld$} of an LD is the set of all permutations of its variables that send  feasible points to feasible points. The symmetry group of the LD relaxation of an ILD is contained in the symmetry group of the ILD. Geyer~\etal~\cite{Geyer2018} provided a method for finding the symmetry group of a linear program (LP).  The symmetry group of an ILD or an 
 LD is related to the symmetry group of an integer linear program (ILP) or an LP as follows.
 If each feasible point of an ILP (LP) is also optimal, then
 the symmetry group of the ILD (LD) of the feasible set of this ILP (LP) coincides with the symmetry group of the ILP (LP).  Hence, the method provided in~\cite{Geyer2018} can be used to find the symmetry group of an LD by applying it to the LP obtained from the LD by making the LD LP's feasible set and the zero function LP's objective function.
  Throughout the paper when we refer to an ILD as an LD we mean the LD relaxation of that ILD. It is shown in Geyer~\etal~\cite{Geyer2018} that  $$S_n\wr S_k\cong G^{\rm iso}(k,n)\leq G^{\rm LD(\ref{ilp:BF})}.$$ 
Moreover, for arbitrary  permutations $h_1,\ldots,h_k$ of the elements of $\{l_1,\ldots,l_n\}$  
and an arbitrary permutation $g$ of the elements of $\{1,\ldots,k\}$, 
 each $((h_1,\ldots,h_k),g)\in G^{\rm iso}(k,n)$ acts transitively on the 
variables of ILD~(\ref{ilp:BF}) by permuting the entries of the frequency vector $\x$ 
according to 
 \begin{equation}\label{eqn:act}
\begin{array}{rll}
((h_1,\ldots,h_k),g)x(i'_1,\ldots,i'_k)&\hspace{-.2cm}=&\hspace{-.2cm}x(i_1,\ldots,i_k),\\
((h_1,\ldots,h_k),g)(i_1,\ldots, i_k)&\hspace{-.2cm}=&\hspace{-.2cm}(i'_1,\ldots, i'_k),
 \end{array}
 \end{equation}
where $(i_1,\ldots,i_k)^{\top}\in \{l_1,\ldots, l_n\}^k$ and $(i'_1,\ldots, i'_k)=(h_1(i_{{g^{-1}}(1)}),\ldots, h_k(i_{{g^{-1}}(k)}))$. Throughout the paper, unless otherwise stated, the action of $S_n\wr S_k$ or  one of its subgroups  on a vector in $\mathbb{C}^{n^k}$ is defined according to equation~(\ref{eqn:act}).

For a subgroup $G$ of the symmetry group of an ILD (ILP), two solutions $\x_1$, $\x_2$ of an ILD (ILP) are called 
{\em isomorphic} with respect to $G$ if there exists some 
$g\in G \, \ni \, g(\x_1)=\x_2$. Margot~\cite{Margot2007}  developed the branch-and-bound with isomorphism pruning algorithm for solving an ILP 
by exploiting a given subgroup $G$ of its symmetry group. An altered version of this algorithm, that finds a set of all non-isomorphic solutions of an ILP with respect to a given subgroup $G$ of its symmetry group, was used 
in~\cite{Bulutoglu2008,Bulutoglu2016} to classify all non-isomorphic $\oand$ for many $k,n,s,\lambda$ combinations. 

For a subgroup $H$ of the symmetry group of an ILD and the constraint $\vv^{\top}\x=c$ for some constant $c\in \mathbb{R}$ implied by the integrality constraints of the ILD, the  non-trivial constraints 
\begin{equation}\label{eqn:HHHH}
(h^{\top}(\vv)-\vv)^{\top}\x=0\quad \text{for}\,\,h\in H
\end{equation}
are valid for the feasible set of the ILD. We call such constraints the 
{\em zero right hand side linear equality constraints associated with $H$.} The valid equalities~(\ref{eqn:HHHH}) based on $H$ put restrictions on the dimension of the convex hull of all feasible solutions of the ILD. It is far from clear what these restrictions would be. We establish  Theorem~\ref{thm:dimfin} that provides such restrictions
for the orthogonal array problem. We also discuss how to find such restrictions as well as the corresponding classes of potentially valid equalities of the form~(\ref{eqn:HHHH}) for a general ILD.
 
 The paper is organized as follows. In Section~\ref{sec:ANOVA}, 
  we review the theory of analysis of variance (ANOVA) by using representation 
  theory~\cite{Diaconis1988}.
 In Section~\ref{sec:J}, we introduce the concept of the $J$-characteristics of an array, and provide a set of necessary and sufficient constraints for an array to be an OA based on its $J$-characteristics. Moreover, we prove that certain constraints  must be satisfied by  the $J$-characteristics of OAs.
 In Section~\ref{sec:dim2}, we determine the decomposition of $\mathbb{Q}^{\X}$ into irreducible subrepresentations under the action of the largest known subgroup of $G^{\text{LD(\ref{ilp:BF}})}$.
  In Section~\ref{sec:dim}, we use representation theory, ANOVA, and the results of Section~\ref{sec:dim2} to  show that the symmetries of $P^{(k,s,\lambda)}_{ n;I}$ drastically decrease the number of possible values of dim($P^{(k,s,\lambda)}_{ n;I}$). By using the $J$-characteristics, we also determine the corresponding sets of linear equality constraints that can potentially  be satisfied by all the points in $P^{(k,s,\lambda)}_{n;I}$. These are the only linear equality constraints up to equivalence 
  that can potentially be implied by the integrality constraints of  ILD~(\ref{ilp:BF}).  Our results imply all  the values of dim($P^{(k,s)}_{n;I})$ in Table~\ref{tab:PnI}.
Moreover, we find additional restrictions on dim($P^{(k,s,\lambda)}_{n;I})$ and larger sets of corresponding linear equality constraints for the $n=2$ and even $s$ cases that possess symmetries that do not necessarily exist
 in the $3\leq n$ or odd $s$ cases. These sets of linear equality constraints are  obtained by taking the union of the sets of linear equality constraints obtained for the general case.
 Again, these are the only linear equality constraints up to equivalence that can potentially be implied by the integrality constraints of  ILD~(\ref{ilp:BF}).
   In Section~\ref{sec:generalization},  we discuss how representation theory can be used to decrease the number of possible values of the dimension of the convex hull of all feasible points of a general ILD with a given subgroup $H$ of its symmetry group.  We also discuss how to generate the corresponding sets of zero right hand side linear equality constraints associated with $H$. These are the only zero right hand side linear equality constraints associated with $H$ up to equivalence that can potentially  be implied by the integrality constraints of the ILD. 
Finally we summarize the main findings of the paper, and emphasize an open problem in representation theory that stems from Theorem~\ref{thm:nointeger} in Section~\ref{sec:generalization}.

Throughout the paper, a lowercase boldfaced letter is a coordinate vector with respect to a fixed basis, and an uppercase boldfaced letter is a matrix. 
The vector $\1_n$ ($\0_n$) is the all $1$s ($0$s) vector of length $n$.
 For a vector $\z \in \mathbb{Q}^n$ and a group $G$ that acts on $\z$ by permuting its entries, $G\z$ is  the orbit of $\z$ under the action of  $G$, that is, $$G\z:=\{ \vv\in\mathbb{Q}^n\ |\ \vv =g(\z) \text{ for some } g\in G \}.$$
For any vector space in the paper, unless otherwise stated, the field of scalars $\mathbb{F}$ is 
$\mathbb{Q},\ \mathbb{R},$ or $\mathbb{C}$. 
Since the complex dot product coincides with the dot product in $\mathbb{R}$ or in $\mathbb{Q}$ we will use 
the complex dot product in place of the dot product for the rest of the paper. Unless otherwise stated, the orthogonal complement $V^{\perp}$
of a vector space $V$ and orthogonal direct sums denoted by ``$\obot$"  are with respect to the complex dot product.
For a set of points $S$ in a vector space over the field of scalars $\mathbb{F}$, 
Span$_{\mathbb F}(S$) is the span, Aff$_{\mathbb F}(S$) is the affine hull, and dim$_{\mathbb F}(S$) is the dimension  of the affine hull  of the vectors in $S$ over $\mathbb{F}$. For a matrix $\A$, Row$_{\mathbb{F}}(\A$) is the row space,  Null$_{\mathbb{F}}(\A$) is the null space, and Col$_{\mathbb{F}}(\A$) is the column space  of $\A$  over ${\mathbb F}$.
 If ${\mathbb F}$ is not provided, then ${\mathbb F}={\mathbb R}$.  
 Finally, Conv($S$) is the convex hull of the points in $S$. 
 \section{The irreducible representations of $\prod_{i=1}^k S_n $ in ANOVA} \label{sec:ANOVA}
We first provide some background material on group representations. When a group $G$ acts on a vector space $V$ over a field $\mathbb{F}$, i.e., there is a homomorphism $\rho:G \rightarrow \text{Aut}_{\mathbb{F}}({V})$ from $G$ into the group of $\mathbb{F}$-linear automorphisms of the vector space $V$, then (by abuse of language) both this homomorphism and $V$ 
under this action are called a {\em representation} of 
$G$~\cite{Diaconis1988, Serre}. The representation $(\rho, V)$ is called  {\em rational}, {\em real}, {\em complex} when $\mathbb{F}$ is
$\mathbb{Q}$,  $\mathbb{R}$,  $\mathbb{C}$, respectively. 
Throughout the paper we will call the representation $(\rho, V)$ an $\mathbb{F}$-representation. 
 A $G$-invariant subspace $W$ of $V$ 
yields by restriction a homomorphism $\rho_{|W}:G\rightarrow  \text{Aut}_{\mathbb{F}}({W})$,  and both $W$ and this homomorphism are called a {\em subrepresentation} of $V$.  

 An $\mathbb{F}$-representation $\rho: G\rightarrow \text{Aut}_{\mathbb{F}}({V})$ 
 is called {\em trivial} if ${\rm dim}(V)=1$ and $\rho (g)$  acts as the identity on $V$ $\forall g\in G$.
 Two $\mathbb{F}$-representations $\rho_1:G\rightarrow  \text{Aut}_{\mathbb{F}}({W})$ and $\rho_2:G\rightarrow  \text{Aut}_{\mathbb{F}}({W'})$ of $G$  
are {\em equivalent} if there is an invertible linear map $\phi:W \rightarrow W'\, \ni \, \phi(\rho_1(g)w)=\rho_2(g)\phi(w)$ $\forall w \in W$ and $g \in G$.
 Clearly, being equivalent is an equivalence relation among all $\mathbb{F}$-representations of a group $G$.   
Representation theory has been developed to find all non-equivalent $\mathbb{F}$-representations of groups. The character $\chi_{\rho}$ of an $\mathbb{F}$-representation 
$\rho: G\rightarrow \text{Aut}_{\mathbb{F}}({V})$ is
defined to be the map $\chi_{\rho} : G \rightarrow \mathbb{F}$ such that 
$\chi_{\rho} (g) = \Tr(\rho(g))$ for $g \in G$, where $\Tr(\rho(g))$ is the trace of the linear transformation $\rho(g)$.
 Two $\mathbb{F}$-representations of a finite group are equivalent if and only if they have the same character~\cite{Fassler}. A $\CCC$-representation $(\rho,V)$ of a finite group $G$ with character $\chi_{\rho}$ is irreducible if and only if $\langle \chi_{\rho} \mid \chi_{\rho} \rangle=(1/|G|)\sum_{g \in G}\overline{\chi_{\rho}(g)}
 \chi_{\rho}(g)=1$, where $\langle \cdot\mid  \cdot \rangle$ is an inner product between functions from $G$ to $\CCC$~\cite{Fassler}.
 
 An $\mathbb{F}$-representation $\rho:G \rightarrow  \text{Aut}_{\mathbb{F}}({V})$ is {\em unitary} with respect to an inner product  $\left\langle\cdot,\cdot\right\rangle$
  defined in $V$
  if $\left\langle\rho(g)\vv,\rho(g)\uuu\right\rangle=\langle\uuu,\vv\rangle$ for all $\uuu,\vv \in V$.
   It is well-known that every $\mathbb{F}$-representation is 
 unitary with respect to some inner product~\cite[Theorem~1 on p.~8]{Diaconis1988}. 
  An $\mathbb{F}$-representation of a group is called a {\em permutation $\mathbb{F}$-representation} if its action on $V$ can be identified 
with permutations of a basis of $V$. 
A permutation $\mathbb{F}$-representation is unitary with respect to the complex dot product.
Let $\mathbb{F}^{\{l_1,\ldots,l_n\}}$ be the set of all 
$\mathbb{F}$-vectors indexed by the symbols $\{l_1,\ldots,l_n\}$ in ILD~(\ref{ilp:BF}). Then 
$\mathbb{F}^n\cong \mathbb{F}^{\{l_1,\ldots,l_n\}}=
\text{Span}_{\mathbb{F}}(e_{l_1},\ldots, e_{l_n})$, where $e_{l_i}$ is the vector indexed by the symbols $\{l_1,\ldots, l_n\}$ such that $e_{l_i}$ is one at the $l_i$th position, and zero elsewhere.
Let $S_{\{{l_1},\ldots, {l_n}\}}$ be the  group of all permutations of $\{{l_1},\ldots, {l_n}\}$. Then  $S_{\{{l_1},\ldots, {l_n}\}}\cong S_n$ acts on the vector space $\mathbb{F}^{\{l_1,\ldots,l_n\}}=
\text{Span}_{\mathbb{F}}(e_{l_1},\ldots, e_{l_n})$ by 
$\pi e_{l_i}= e_{\pi(l_i)}$ for each $\pi \in  S_{\{{l_1},\ldots, {l_n}\}}$. The action of the group $S_{\{{l_1},\ldots, {l_n}\}}\cong S_n$ is a permutation $\mathbb{F}$-representation of $S_{\{{l_1},\ldots, {l_n}\}}$, and
the subspace ${\rm Span}_{\mathbb{F}}(\1_n)$ is the trivial $\mathbb{F}$-representation of $S_{\{{l_1},\ldots, {l_n}\}}$ appearing as a subrepresentation.  
If an $n$-dimensional  $\mathbb{F}$-representation of a group $\rho:G\rightarrow  \text{Aut}_{\mathbb{F}}({V})$ cannot be further decomposed into invariant subspaces by employing a change of bases, i.e., there exists no invariant subspaces $V_1\neq \{\0_n\}$ and $V_2\neq\{\0_n\}$ of $V$ such that 
 $V=V_1 \oplus V_2$ and 
$\rho_{|V_i}:G\rightarrow  \text{Aut}_{\mathbb{F}}({V_i})$
for $i=1,2$ are both $\mathbb{F}$-representations of $G$, then 
$\rho:G\rightarrow  \text{Aut}_{\mathbb{F}}({V})$ is called an {\em irreducible $\mathbb{F}$-representation} of $G$.
The $n-1$-dimensional subspace $\1_n^{\perp}$ is an irreducible  $\mathbb{C}$-representation of $S_n$~\cite{Diaconis1988}, and consequently $\1_n^{\perp}$ is an irreducible  $\mathbb{F}$-representation of $S_n$ for any $\mathbb{F}$ such that $\mathbb{Q}\subseteq \mathbb{F}\subseteq \mathbb{C}$.

For the rest of the paper, let $\X$ be the $n^k \times k$ array, where the  rows of  $\X$ consist of each of the distinct $n^k$ symbol combinations from $\{l_1,\ldots,l_n\}^k$ ordered lexicographically.
Let $\mathbb{Q}^{\X}$ ($\mathbb{C}^{\X}$) be the vector space of all functions from $\{\text{rows of }\X\}$ to ${\mathbb Q}$ (${\mathbb C}$). Then $$\mathbb{Q}^{\X}\cong(\mathbb{Q}^{n})^{\otimes k}, \quad \mathbb{C}^{\X}\cong(\mathbb{C}^{n})^{\otimes k},$$
$\mathbb{Q}^{\X}=\text{Span}_{\mathbb{Q}}(e_{\x_1},\ldots, e_{\x_{n^k}})$, and $\mathbb{C}^{\X}=\text{Span}_{\mathbb{C}}(e_{\x_1},\ldots, e_{\x_{n^k}})$,
where $\x_i$ is the $i$th row $\X$, and $e_{\x_i}\in \mathbb{Q}^{\X}$ is the function that takes the value $1$ at $\x_i$, and zero at every $\x_j\neq \x_i$ such that $\x_j$ is a row of $\X$. Let $S_{\{l_1,\ldots,l_n\}_j} $ be the group of all permutations of the symbols on the $j$th column of $\X$.
Then the group $\prod_{i=1}^k S_{\{l_1,\ldots,l_n\}_i}\cong \prod_{i=1}^k S_n$ acts  on the elements of $\{e_{\x_1},\ldots, e_{\x_{n^k}}\}$ by acting on the columns 
of $\X$,  and the 
resulting action of $\prod_{i=1}^k S_{\{l_1,\ldots,l_n\}_i}$ on $\mathbb{Q}^{\X}$ and $\mathbb{C}^{\X}$ are both permutation representations.

ANOVA is a decomposition of 
$\mathbb{Q}^{\X}\cong(\mathbb{Q}^n)^{\otimes k}$ ($\mathbb{C}^{\X}\cong(\mathbb{C}^n)^{\otimes k}$)
 into $2^k$ mutually orthogonal subspaces~\cite{Takemura1983}. 
These subspaces can be found by first considering the case $k=1$.
For $k=1$, $\mathbb{Q}^{\X}\cong \mathbb{Q}^n$ decomposes into the direct sum of two subspaces that are invariant under the action of
 $ S_{\{l_1,\ldots,l_n\}}\cong S_n$, i.e.,
$${\mathbb Q}^n\cong\mathbb{Q}^{\X}={\rm Span}_{\mathbb{Q}}(\1_n) \obot \1_n^{\perp},$$
where $S_{\{l_1,\ldots,l_n\}}\cong S_n$ permutes the symbols $\{l_1,\ldots,l_n\}$ in the column of $\X$. 
 For $k=2$, and $i \in \{1,2\}$ let $S_{\{l_1,\ldots, l_n\}_i}$ 
 be the group of all permutations of the symbols 
 $\{l_1,\ldots,l_n\}$ in the $i$th column of $\X$. Then we get the following orthogonal decomposition into irreducible invariant subspaces under the action of $S_{\{l_1,\ldots,l_n\}_1} \times S_{\{l_1,\ldots,l_n\}_2}$ as 
 in~\cite[p.~155]{Diaconis1988}, 
{\footnotesize
\begin{equation*} \label{eqn:twofactor}
\begin{array}{cccccc}
 {\mathbb Q}^n\otimes {\mathbb Q}^n&\hspace{-.35cm} \cong  \mathbb{Q}^{\X} &\hspace{-.35cm}\cong({\rm Span}_{\mathbb{Q}}{(\1_n)}_1\otimes {\rm Span}_{\mathbb{Q}}{(\1_n)}_2)& \hspace{-.35cm} \obot ({\rm Span}_{\mathbb{Q}}{(\1_n)}_1\otimes {(\1_n^{\perp})}_2) & \hspace{-.35cm}  \obot ({(\1_n^{\perp})}_1 \otimes {\rm Span}_{\mathbb{Q}}(\1_n)_2) &  \hspace{-.35cm} \obot  ({(\1_n^{\perp})}_1 \otimes {(\1_n^{\perp})}_2),\\
 & n^2 & 1 & n-1 & n-1 &(n-1)^2,
  \end{array}
\end{equation*}}where the values below each subspace is its dimension.
By using tensor powers, and taking into account the multiplicities of 
each non-equivalent irreducible invariant subspace that appears in this decomposition, we get 
\begin{equation} \label{eqn:twofactor2}
  ({\mathbb Q}^n)^{\otimes 2}\cong \mathbb{Q}^{\X} \cong ({\rm Span}_{\mathbb{Q}}(\1_n))^{\otimes 2}  \obot  (\1_n^{\perp} \otimes {\rm Span}_{\mathbb{Q}}(\1_n)) \obot ({\rm Span}_{\mathbb{Q}}(\1_n)\otimes \1_n^{\perp}) \obot ( \1_n^{\perp})^{\otimes 2}.
\end{equation} 
To generalize this result, we need the following lemma 
 from~\cite{Diaconis1988}. 
\begin{lem}\label{lem:tensor}
Let $G_1$ and $G_2$ be finite groups. Let $\rho_1:G_1\rightarrow {\rm GL}(V_1)$ and $\rho_2:G_2\rightarrow {\rm GL}(V_2)$ 
be $\mathbb{C}$-representations. 
Then, for the $\mathbb{C}$-representation $\rho_1 \otimes \rho_2:G_1\times G_2\rightarrow {\rm GL}(V_1\otimes V_2)$ defined by 
$$ \rho_1\otimes\rho_2(s,t)(\vv_1\otimes\vv_2)=\rho_1(s)(\vv_1)\otimes\rho_2(t)(\vv_2),$$ the following hold. 
\begin{enumerate}
\item If $\rho_1$ and $\rho_2$ are irreducible, then $\rho_1 \otimes \rho_2$ is irreducible.
\item Each irreducible $\mathbb{C}$-representation of $G_1\times G_2$ is equivalent to a $\mathbb{C}$-representation  $\rho_1 \otimes \rho_2$, where for $i=1,2$ $\rho_i$ is an irreducible $\mathbb{C}$-representation of $G_i$.
\end{enumerate}
\end{lem}

Throughout the rest of the paper let $[k] = \{1,\ldots,k\}$. Now, we prove the following theorem.
\begin{thm}\label{thm:generalk}
Let  $u\subseteq [k]$, $U_{0,i} =
 {\rm Span}_{\mathbb{Q}}(\1_n )$, $U'_{0,i} =
 {\rm Span}_{\mathbb{C}}(\1_n )$, $U_{1,i} = \1_n^\perp \subset \mathbb{Q}^n$, and  $U'_{1,i} = \1_n^\perp \subset \mathbb{C}^n$ for $i\in [k]$. For $u \subseteq [k]$, let
$L_u := U_{\varepsilon_1,1} \otimes U_{\varepsilon_2,2} \otimes \cdots \otimes U_{\varepsilon_k,k}$ and $L'_u := U'_{\varepsilon_1,1} \otimes U'_{\varepsilon_2,2} \otimes \cdots \otimes U'_{\varepsilon_k,k}$ with $\varepsilon_i = 1$ when $i \in u$, and $\varepsilon_i = 0$ otherwise. Then for general $k$,  orthogonal decomposition~(\ref{eqn:twofactor2}) of $\mathbb{Q}^{\X}$ and $\mathbb{C}^{\X}$ into irreducible invariant subspaces under the action of $\prod_{i=1}^k S_{\{l_1,\ldots,l_n\}_i}$ are
\begin{equation} \label{eqn:all}
 ({\mathbb Q}^n)^{\otimes k}\cong\mathbb{Q}^{\X}  = \bigotimes_{i=1}^k
 ({\rm Span}_{\mathbb{Q}}(\1_n)
\obot \1_n^{\perp})_i \cong
\bigobot_{u \subseteq [k]} L_u
\end{equation} 
and
\begin{equation} \label{eqn:all2}
 ({\mathbb C^n})^{\otimes k}\cong\mathbb{C}^{\X}  
 = \bigotimes_{i=1}^k({\rm Span}_{\mathbb{C}}(\1_n)
\obot\1_n^{\perp})_i \cong  \bigobot_{u \subseteq [k]} L'_u,
\end{equation} 
where the mutual orthogonality of subspaces $L_u \ni u \subseteq [k]$ and $L'_u \ni u \subseteq [k]$ is with respect to  the complex dot product.
\end{thm}
\begin{pf}
First, equations~(\ref{eqn:all}) and~(\ref{eqn:all2}) are clear by the properties of tensor products and direct sums of vector spaces. Also, observe that $S_{\{l_1,\ldots,l_n\}_{i}}\cong  S_n$ for all possible $i$.
For the $i$th column of $\X$, let $\rho_{0,i}$ and $\rho_{1,i}$ be such that
$$\rho_{0,i}\obot\rho_{1,i}: S_{\{l_1,\ldots,l_n\}_i}\rightarrow 
\text{GL}(({\rm Span}_{\mathbb{F}}(\1_n) 
\obot \1_n^{\perp})_i),$$ and $\rho_{0,i}$ and $\rho_{1,i}$ are the irreducible $\mathbb{F}$-representations of $S_{\{l_1,\ldots,l_n\}_i}$ corresponding to
 ${\rm Span}_{\mathbb{F}}(\1_n)$ and 
$\1_n^{\perp}$ in $({\rm Span}_{\mathbb{F}}(\1_n) 
\obot \1_n^{\perp})_i$, where 
$\mathbb{F}=\mathbb{Q}$ or $\mathbb{F}=\mathbb{C}$. 

 For $\mathbb{F}=\mathbb{C}$ let
$$\rho_u := \rho_{\varepsilon_1,1} \otimes \rho_{\varepsilon_2,2} \otimes \cdots \otimes \rho_{\varepsilon_k,k}:\prod_{i=1}^k S_{\{l_1,\ldots,l_n\}_i}\rightarrow 
\text{GL}\left( L'_u \right).$$ Then by the properties of tensor products and direct sums of representations
$$\bigotimes_{i=1}^k\left(\rho_{0,i}\obot\rho_{1,i}\right)=\bigobot_{u \subseteq [k]} \rho_u:\prod_{i=1}^k S_{\{l_1,\ldots,l_n\}_i}\rightarrow \text{GL}\left(\bigotimes_{i=1}^k({\rm Span}_{\mathbb{C}}(\1_n)
\obot \1_n^{\perp})_i\right)= \text{GL}\left(\bigobot_{u \subseteq [k]} L'_u \right).$$
Moreover, 
by using induction on $k$, and applying  Lemma~\ref{lem:tensor} $2^{k}$ times, we get that each $L'_u$
 is an irreducible $\mathbb{C}$-representation of 
$\prod_{i=1}^k S_{\{l_1,\ldots,l_n\}_i}\cong \prod_{i=1}^k S_{n} $. The proof for $\mathbb{Q}^{\X}$ is obtained by first replacing $\mathbb{C}$ with $\mathbb{Q}$ and $L'_u$ with $L_u$ in the first part of the proof, and observing that $L_u$ is invariant under the action of $\prod_{i=1}^k S_{\{l_1,\ldots,l_n\}_i}$, where the irreducibility of $L_u$ follows from the irreducibility of $L'_u$ and the fact that $\mathbb{Q} \subset \mathbb{C}$.\qed
\end{pf}
Decomposition~(\ref{eqn:all})~((\ref{eqn:all2})) is known as the {\em ANOVA decomposition of} $({\mathbb Q}^n)^{\otimes k}$ $(({\mathbb C^n})^{\otimes k})$~\cite{Takemura1983}.
Using a basis that allows  decomposition~(\ref{eqn:all})~((\ref{eqn:all2})) to express a function $f((i_1,\ldots,i_k)) \in \mathbb{Q}^{\X}$ ($f((i_1,\ldots,i_k)) \in \mathbb{C}^{\X}$) is called an {\em ANOVA decomposition of} $f((i_1,\ldots,i_k))$. 
 The generalization of the ANOVA decomposition of $({\mathbb Q}^n)^{\otimes k}$ ($({\mathbb C}^n)^{\otimes k}$) to the ANOVA decomposition of 
 $\otimes_{i=1}^k{\mathbb Q}^{n_i}$ ($\otimes_{i=1}^k{\mathbb C}^{n_i}$) is straightforward~\cite{Takemura1983}, and each of the $2^k$ subspaces that appear in this decomposition is equivalent to an irreducible $\mathbb{Q}$-representation ($\mathbb{C}$-representation) of $\prod_{i=1}^k S_{n_i}$~\cite{Diaconis1988}.
\section{J-characteristics}\label{sec:J}
An array $\D$ of  $N$ rows and $k$ columns with entries from the set 
$\{l_1,\ldots,l_n\}$  is called an 
{\em $N$ row, $k$ column, $n$-symbol array}.  
For a given  $\D$, let $x(i_1, \ldots, i_k)$ be the number of times the  symbol combination $(i_1,\ldots,i_k)$ such that  
$(i_1,\ldots,i_k)^{\top}\in \{l_1,\ldots,l_n\}^k$  appears in $\D$. 
 Then the grand mean based on $\D$ 
 is defined to be
\begin{equation}\label{eqn:grandmean}
x_{\emptyset}(i_1, \ldots, i_k) = n^{-k} 
\sum_{i_1, \ldots, i_k} x(i_1, \ldots, i_k)=\frac{N}{n^k},
\end{equation}
 and for $u\subseteq [k]$ the interaction $x_u (i_1, \ldots, i_k)$ involving the columns indexed by the indices in $u$ is defined by
\begin{equation}\label{eqn:unnumbered}
x_u (i_1, \ldots, i_k) = 
n^{-k+|u|}
\sum_{ \{i_j \,\mid \, j \not\in u\}}  
x(i_1, \ldots, i_k) 
- \sum_{ v \subsetneq u } x_v (i_1, \ldots, i_k).
\end{equation}
Moreover, \begin{equation}\label{eqn:ANOVA}
	x(i_1, \ldots, i_k) = \sum_{ u \subseteq [k] } x_u (i_1, \ldots, i_k),
\end{equation}
and equation~(\ref{eqn:ANOVA}) is  the  ANOVA decomposition of $x(i_1, \ldots, i_k)$~\cite{Christensen}. 
The $J$-characteristics in~\cite[p.63]{Lekivetz2011} are defined as 
\begin{equation}\label{eqn:Jr}
J_u^{\x} (i_1, \ldots, i_k) = n^k x_u (i_1, \ldots, i_k),
\end{equation}
where $\x \in \mathbb{Z}^{n^k}$ is indexed by the elements in $\{l_1,\ldots,l_n\}^k$
whose  $(i_1, \ldots, i_k)$th entry is  $x(i_1, \ldots, i_k)$.
 By induction, 
 each of $x_u (i_1, \ldots, i_k)$ and $J_u^{\x} (i_1, \ldots, i_k)$
is a function of the indices indexed by the elements in $u$ only, and does not depend on the indices indexed by the elements in $[k]\backslash u$.
By equations~(\ref{eqn:ANOVA})~and~(\ref{eqn:Jr}), we have
\begin{equation}\label{eqn:Ju}
n^k x(i_1, \ldots, i_k) 
= \sum_{ u \subseteq [k] } J^{\x}_u (i_1, \ldots, i_k).
\end{equation}
The following lemma shows that
$J^{\x}_u (i_1, \ldots, i_k)$ for each $u \subseteq [k]$ depends only on  
$\D_u$, where 
$\D_u$  is obtained from $\D$ by deleting the columns of $\D$ indexed by the indices in $[k]\backslash u$. 
 \begin{lem}\label{lem:projD}
 Let $\D,\,\D'$ be  $n$-symbol arrays with $k,\,k'$ columns such that $\D'$ is not necessarily equal to $\D$. For each symbol combination $(i'_1,\ldots,i'_{k'})^{\top} \in
  \{l_1,\ldots,l_n\}^{k'}$, let $x'(i'_1, \ldots, i'_{k'})$  be the number of times 
 $(i'_1,\ldots,i'_{k'})$ 
 appears as a row of $\D'$. Let $\x$ be as in equation~(\ref{eqn:grandmean}), and  $\x' \in \mathbb{Z}^{n^{k'}}$ be indexed by the elements in $\{l_1,\ldots,l_n\}^{k'}$, 
where  $(i_1, \ldots, i_{k'})$th entry of $\x'$ is  $x'(i_1, \ldots, i_{k'})$.
 Let $u=\{j_1,\ldots, j_{|u|}\}\subseteq [k]$, 
 $u'=\{j'_1,\ldots, j'_{|u'|}\}\subseteq [k']$ be such that $$\text{the multiset of rows of } \D_u=\text{the multiset of rows of }\D'_{u'},$$  then $J^{\x}_u (i_1, \ldots, i_k)=J^{\x'}_{u'} (i'_1, \ldots, i'_{k'})$.
 \end{lem}
 \begin{pf}
 The proof follows by induction on $|u|=|u'|$. \qed
 \end{pf}
 The concept of $J$-characteristics can also be described by using the $k$-way layout fixed effects interpolation model in statistics for an all possible combinations experiment with $k$ columns, each column having $n$ distinct symbols from $\{l_1,\ldots,l_n\}$ replicated $m=1$ times, i.e., each of the $n^k$ symbol combinations appearing exactly $m$ times for $m=1$. In particular, the $3$-way layout fixed effects  model for the response variable $Y_{{i_1}{i_2}{i_3}j}$ of such an experiment for general 
 $m$ has the form 
 \begin{equation}\label{eqn:3way}
Y_{{i_1}{i_2}{i_3}j}=\alpha^{\emptyset}+\alpha^1_{{i_1}}+\alpha^2_{{i_2}}+\alpha^3_{{i_3}}+
\alpha^{12}_{{i_1}{i_2}}+\alpha^{13}_{{i_1}{i_3}}+\alpha^{23}_{{i_2}{i_3}}+\alpha^{123}_{{i_1}{i_2}{i_3}}+
\epsilon_{{i_1}{i_2}{i_3}j}
 \end{equation} for $(i_1,i_2,i_3,j)\in \{l_1,\ldots,l_n\}^3\times [m]$, 
where $\epsilon_{{i_1}{i_2}{i_3}j}$ are identically independently distributed as
 N$(0,\sigma^2)$ for some $\sigma^2\geq 0$, and the following equations
 \begin{eqnarray}\label{eqn:side}
 \begin{array}{lllll}
& \sum_{i_1\in\{\l_1,\ldots,l_n\}}\alpha_{{i_1}}^1&\hspace{-.3cm}=0,
&\sum_{i_2\in \{\l_1,\ldots,l_n\}}\alpha_{{i_2}}^2&\hspace{-.3cm}=0, \\ 
& \sum_{i_3\in\{\l_1,\ldots,l_n\}}\alpha_{{i_3}}^3&\hspace{-.3cm}=0, 
\quad &\sum_{i_2\in \{\l_1,\ldots,l_n\}}
\alpha_{{i_1}{i_2}}^{12}&\hspace{-.3cm}=0\text{ for each ${i_1}$,}\\
&\sum_{i_1\in\{\l_1,\ldots,l_n\}}
 \alpha_{{i_1}{i_2}}^{12}&\hspace{-.3cm}=0\text{ for each ${i_2}$,}
 \quad &\sum_{i_3\in \{\l_1,\ldots,l_n\}}
\alpha_{{i_1}{i_3}}^{13}&\hspace{-.3cm}=0\text{ for each ${i_1}$,}\\
&\sum_{i_1\in \{\l_1,\ldots,l_n\}}
\alpha_{{i_1}{i_3}}^{13}&\hspace{-.3cm}=0\text{ for each ${i_3}$,}\quad
 &\sum_{i_2\in \{\l_1,\ldots,l_n\}}
\alpha_{{i_2}{i_3}}^{23}&\hspace{-.3cm}=0\text{ for each ${i_3}$,}\\
&\sum_{i_3\in \{\l_1,\ldots,l_n\}}
\alpha_{{i_2}{i_3}}^{23}&\hspace{-.3cm}=0\text{ for each ${i_2}$,}\quad
&\sum_{i_3\in \{\l_1,\ldots,l_n\}}
\alpha_{{i_1}{i_2}{i_3}}^{123}&\hspace{-.3cm}=0\text{ for each $({i_1},{i_2})$ tuple,}\\
&\sum_{i_2\in \{\l_1,\ldots,l_n\}}
\alpha_{{i_1}{i_2}{i_3}}^{123}&\hspace{-.3cm}=0\text{ for each $({i_1},{i_3})$ tuple,}
\quad &\sum_{i_1\in \{\l_1,\ldots,l_n\}}
\alpha_{{i_1}{i_2}{i_3}}^{123}&\hspace{-.3cm}=0\text{ for each $({i_2},{i_3})$ tuple,}
\end{array}
 \end{eqnarray}
 are satisfied by the {\em main effect parameters} (parameters with a single index) and {\em interaction parameters} (parameters with more than one index) of the model.
 Equations~(\ref{eqn:side}) are called the {\em side constraints}.
Generalization to $k$-way layout is straightforward, and in this case, the ${ r\choose r-1 }$ side constraints 
 for $\alpha_{{i_1}\ldots {i_r}}^{1\ldots r}$ are the same as the equality constraints in ILD~(\ref{ilp:BF}) for an OA$(N,r,n,r-1)$ except the right hand side vector for the equality constraints is $\0_q$ instead of $N/n^{r-1}\1_q$, where $q={r \choose r-1} n^{r-1}$. 
 Given the observed values $y_{{i_1}{i_2}{i_3}j}$ of $Y_{{i_1}{i_2}{i_3}j}$, the ordinary least squares method for the fixed effects model seeks to find estimates for the main effect and interaction  parameters  by solving
 \begin{equation}\label{eqn:optcon}
 \begin{array}{rl}
& \min \quad \sum_{{i_1},{i_2},{i_3},j}  (y_{{i_1}{i_2}{i_3}j}-\alpha^{\emptyset}-\alpha^1_{{i_1}}-\alpha^2_{{i_2}}-\alpha^3_{{i_3}}-
\alpha^{12}_{{i_1}{i_2}}-\alpha^{13}_{{i_1}{i_3}}-\alpha^{23}_{{i_2}{i_3}}-\alpha^{123}_{{i_1}{i_2}{i_3}})^2  \\
& \mbox{s.t.: \quad equations~(\ref{eqn:side}) are satisfied.}
\end{array}
\end{equation}
Optimization problem~(\ref{eqn:optcon}) is convex, and has a unique solution attaining the global minimum. This solution provides the estimates for the main effects and interaction parameters in  model~(\ref{eqn:3way}).
  In fact, for $u=\{j_1,\ldots,j_{|u|}\}\subseteq [k]$, the $n^{|u|}$ parameter estimates for the main effect and interaction parameters involving the columns indexed by the elements in $u$ in the  $k$-way layout 
fixed effects model for $x(i_1,\ldots,i_k)$ in~(\ref{eqn:ANOVA}) are $$x_u (i_1, \ldots, i_k)=
\frac{J_u^{\x} (i_1, \ldots, i_k)}{n^k},$$see~\cite{Christensen}. 

Geyer~\etal~\cite{Geyer2018} used a different definition of the $J$-characteristics for arrays with symbols from $\{-1,1\}$.
Next, we provide a simplification of the $J$-characteristics 
in~\cite{Lekivetz2011} for such arrays. This simplification will be used to prove that the definition of the $J$-characteristics used in~\cite{Geyer2018} is consistent with that in~\cite{Lekivetz2011}. 
However, we first need the following lemma obtained by setting $v=2$, and
 replacing 
$\{1,2 \}$ with $\{-1, 1\}$, $t$ with $s$,  $\x$ with $\vv$, and $\y$ with $\w$ in Lemma 2 
of~\cite{Rosenberg1995}.
\begin{lem}\label{lem:Ros}
Let $\{a_c\}$ be such that
\[
a_0=\lambda, \quad a_c=\lambda-\sum_{e=0}^{c-1} a_e {k-s \choose c-e} \quad \text{for } c \geq 1. 
\]
Let $\z$, $\vv$, and $\w$ be row vectors such that  $\z^{\top}$, $\vv^{\top}$, and $\w^{\top} \in \{-1,1\}^k$ with $0 \leq d(\z,\vv)\leq s$, where 
$d(\z,\vv)$ is the number of non-zero entries in $\z-\vv$, i.e., the Hamming distance between $\z$ and $\vv$. 
For a fixed $\z$, let $I_{\vv}=\{i \in [k]:v_i \neq z_i\}$ and $J_{\vv}=\{ \w^{\top} \in  \{-1,1\}^k:w_i=v_i \ \ \forall i \in I_{\vv} \}$.
 Then
\begin{equation*} \label{eqn:Ros}
\begin{aligned}
&N_{\vv}=a_{s-d(\z,\vv)} +(-1)^{s-d(\z,\vv)+1}\sum_{\w\in J_{\vv} \atop 
d(\z,\w)>s} { d(\z,\w)-d(\z,\vv)-1 \choose s-d(\z,\vv)} N_{\w},\\ 
&N_{\w}\geq0, \text{ for $\w\, \ni \, d(\z,\w)>s$},
\end{aligned}
\end{equation*}
where $N_{\vv}$, $N_{\w}$ are the number of times the symbol combinations $\vv$, $\w$  appear in a hypothetical $\oal$.
\end{lem}
The following lemma provides a simplification of the $J$-characteristics $J_u^{\x}(i_1,\ldots,i_k)$ for $2$-symbol arrays with 
symbols from $\{-1,1\}$. 
\begin{lem}\label{lem:RosJ}
For a given $N$ row, $k$ column array $\D$, let $\x \in \mathbb{Z}^{2^k}$ be such that $x(i_1, \ldots, i_k)$  is the number of times the  symbol combination $(i_1,\ldots,i_k)$ with $(i_1,\ldots,i_k)^{\top} \in \{-1,1\}^k$ appears as a row of $\D$.
 For each $u= \{j_1,\ldots,j_{|u|} \} \subseteq [k]$, let $(i_1,\ldots,i_k)_u=(i_{j_1},\ldots,i_{j_{|u|}})$. Then, 
\begin{equation} \label{eqn:RosJ}
\begin{aligned}
J_u^{\x}(i_1,\ldots,i_k)=(-1)^{|u|-d\left(-\1_{|u|}^{\top},(i_1,\ldots,i_k)_u\right)}
J_u^{\x}(1,\ldots,1).
\end{aligned}
\end{equation}
\end{lem}
\begin{pf}
 Let $(i_1,\ldots,i_k)_u=(i_{j_1},\ldots,i_{j_{|u|}})$. Then, $J_u^{\x}(i_1,\ldots,i_k)$ is a function of $(i_{j_1},\ldots,i_{j_{|u|}})$, and there are $2^{|u|}$ distinct assignments for the values of $J_u^{\x}(i_1,\ldots,i_k)$. 
 Moreover, the main effect parameter estimates if $|u|=1$, and the interaction parameter estimates involving the columns indexed by the elements in $u$ if $|u|>1$
 $$\frac{J_u^{\x}(i_1,\ldots,i_k)}{2^k}$$ 
  in the  $k$-way layout fixed effects
  model for $x(i_1,\ldots,i_k)$ 
   must satisfy the side constraints, i.e.,
   the equality constraints in ILD~(\ref{ilp:BF}), with 
 $s=|u|-1$, $\lambda=0$, $n=2$, and $k=|u|$.
 Then, $2^kJ_u^{\x}(i_1,\ldots,i_k)/2^k=J_u^{\x}(i_1,\ldots,i_k)$ must also satisfy the same constraints as the right hand side of each of these constraints is $0$.
 Hence, the result follows  from Lemma~\ref{lem:Ros} by taking $\z=
   -\1_{|u|}^{\top}$, $a_c=\lambda=0$ for $c\geq 0$, $s=|u|-1$, and $k=|u|$. \qed
\end{pf}
 The following definition of the $J$-characteristics was used
  in~\cite{Geyer2018}.  
 \begin{defn}\label{def:Geyer}
Let $\D=(d_{ij})$ be an $N$ row, $k$ column array with symbols from 
$\{-1,1\}$. Let $r \in [k]$ and $u=\{j_1,\ldots,j_r\}\subseteq [k]$.   Then the integers
\begin{equation*}
J_r(u)(\D):=\sum_{i=1}^N{\prod_{j \in \ell}{d_{ij}}}
\end{equation*}
are called the {\it $J$-characteristics} of $\D$. (For $r=0$, $J_0(\emptyset)(\D):=N$.)
\end{defn} 
 Let the column vectors of $\ZZ^{\top}=[\z_1\,  \cdots\, \z_k]^{\top}$ be all $2^k$  vectors in $\{-1,1\}^k$, where $\ZZ$ is constructed the way $\CC$ is constructed in~\cite{Stufken2007}. For distinct   
 $\{j_1, \ldots , j_r\} \subseteq [k]$ with $r\geq2$, let  $\mathbf{z}_{j_1, \ldots, j_r}$  be the 
 {\em $r$-way Hadamard product} $\mathbf{z}_{j_1} \odot \cdots \odot \mathbf{z}_{j_r}$, where for $p \in [2^k]$ the $p$th row of the vector $\mathbf{z}_{j_1} \odot \cdots \odot \mathbf{z}_{j_r}\in \{-1,1\}^{2^k}$ is
the product of the entries on the $p$th row of the matrix   $[\mathbf{z}_{j_1}\,  \cdots\,  \mathbf{z}_{j_r}]$.
Let $\x \in \mathbb{C}^{\ZZ}$ and 
   $\mathbf{H}$  be the  $2^k \times 2^k$  matrix
\begin{equation}\label{eqn:M}
\mathbf{H} = \begin{bmatrix}
\1^{\top}_{2^k}\\
\mathbf{z}_1^{\top} \\
\vdots \\
\mathbf{z}_k^{\top} \\
\mathbf{z}_{1, 2}^{\top} \\
\vdots \\
\mathbf{z}_{1, \dots, k}^{\top} \end{bmatrix}.
\end{equation}
Then the rows of $\HH$ are orthogonal~\cite{Stufken2007}. Consequently 
$\HH^{\top}\HH=\HH\HH^{\top}=2^k\I$, and $\HH^{-1}=(1/2^k)\HH^{\top}$.  
Let $\x$ be such that  $x_p$ for $p \in [2^k]$ is 
the number of times the $p$th row of $\ZZ$ appears in the $N \times k$ array $\D$ with symbols from 
$\{-1,1\}$. 
Define
 \begin{equation}\label{eqn:previous}
 \J^{\x}=(J_0^{\x}(\emptyset), J_1^{\x}(\{1\}),\ldots,J_1^{\x}(\{k\}),J_2^{\x}(\{1,2\}),\ldots,J_k^{\x}(\{1,\ldots,k\}))^{\top}
\end{equation}
via \begin{equation}\label{eqn:H}
\J^{\x}=\mathbf{H}\x.
\end{equation}
Then the entries of $\J^{\x}$
are the corresponding $J$-characteristics of $\D$. 
 By multiplying both sides of equation~(\ref{eqn:H}) by $(1/2^k) \mathbf{H}^{\top}$
 we get 
 \begin{equation} \label{eqn:xfromj}
 \x=\frac{1}{2^k}\mathbf{H}^{\top}\J^{\x}.
 \end{equation}

We next prove that Definition~\ref{def:Geyer} is consistent with the
definition of the $J$-characteristics in~\cite{Lekivetz2011}. 

\begin{lem}\label{lem:consistent}
Let $\D=(d_{ij})$ be an $N$ row, $k$ column array with symbols from $\{-1,1\}$. Let $r \in [k]$ and $u=\{j_1,\ldots,j_r\}\subseteq [k]$.
Let $\x \in \mathbb{Z}^{2^k}$ be such that $x(i_1, \ldots, i_k)$  is the number of times the  symbol combination $(i_1,\ldots,i_k)$ with $(i_1,\ldots,i_k)^{\top} \in \{-1,1\}^k$ appears as a row of $\D$.
Then $$J_r(u)(\D)=J_{u}^{\x}(1,\ldots,1).$$
 \end{lem}
\begin{pf}
Let $\J^{\x}$ be as in equation~(\ref{eqn:previous}). Then by equation~(\ref{eqn:xfromj})
\begin{equation}\label{eqn:JJ}
2^k\x=\mathbf{H}^{\top}\J^{\x}.
\end{equation}
Moreover, by equations~(\ref{eqn:ANOVA}) and~(\ref{eqn:RosJ}) 

\begin{align}\label{eqn:twotothekx}
\begin{split}
2^kx(i_1,\ldots, i_k)=\sum_{u \subseteq [k]}2^k x_u(i_1,\ldots,i_k)&=
\sum_{u \subseteq [k]}J_u^{\x}(i_1,\ldots,i_k)\\ &=\sum_{u \subseteq 
[k]}(-1)^{|u|-d\left(-\1_{|u|}^{\top},(i_1,\ldots,i_k)_u\right)}
J_u^{\x}(1,\ldots,1).
\end{split}
\end{align}
Let
\begin{equation}\label{eqn:previoushat}
 \hat{\J}^{\x}=(J_{\emptyset}^{\x}(\1_k^{\top}), J_{\{1\}}^{\x}(\1_k^{\top}),\ldots,
 J_{\{k\}}^{\x}(\1_k^{\top}),J_{\{1,2\}}^{\x}(\1_k^{\top}),
 \ldots,J_{\{1,\ldots,k\}}^{\x}(\1_k^{\top}))^{\top}.
\end{equation}
Now, equations~(\ref{eqn:twotothekx}) and~(\ref{eqn:previoushat})
imply 
\begin{equation}\label{eqn:Jhat}
2^k\x=\mathbf{H}^{\top}\hat{\J}^{\x}.
\end{equation}
Then by equations~(\ref{eqn:JJ}) and~(\ref{eqn:Jhat})
\begin{equation*}\label{eqn:finalJhat}
2^k\x=\mathbf{H}^{\top}\J^{\x}=\mathbf{H}^{\top}\hat{\J}^{\x}
 \quad \Rightarrow \quad \J^{\x}=\hat{\J}^{\x}.  \qqed
\end{equation*}
\end{pf}
The following lemma from~\cite[p.67]{Lekivetz2011} follows from the properties of OAs and the
fact that the $J$-characteristics of an array $\D$ are  its 
 coordinates   with respect to an orthogonal basis that allows the ANOVA  decomposition~(\ref{eqn:all}). 
\begin{lem}\label{lem:ii}
Let $\D$ be an $N$ row, $k$ column array with entries from $\{l_1,\ldots,l_n\}$, then the following hold.
\\
(i)  $\D$ is uniquely determined by its $J$-characteristics up to permutations of its rows, and vice versa. \\
(ii) $\D$ is an OA of strength $s$ if and only if
$J^{\x}_u =0$ $\forall u\subseteq [k] \, \ni \, 1\leq|u| \leq s$. \\
\end{lem}
First, we prove two combinatorial identities needed for the next theorem.
\begin{lem}\label{lem:small}
Let $k$ and $s$ be positive integers such that $r=k-s\geq 2$. Then
$$
\sum_{i=0}^{r-1} (-1)^{i+1}{s+i \choose i}{s+r-1 \choose s+i}=0.
$$
\end{lem}
\begin{pf}
\begin{align*}
 \sum_{i=0}^{r-1} (-1)^{i+1}{s+i \choose i}{s+r-1 \choose s+i}&=
\\
 \frac{(s+r-1)(s+r-2)\cdots (r)}{s!}\sum_{i=0}^{r-1} (-1)^{i-1} \frac{(r-1)!}{ i! (r-1-i)!}&=\\ 
 \frac{-(s+r-1)(s+r-2)\cdots (r)}{s!}
\sum_{i=0}^{r-1}(-1)^i{r-1 \choose i}&=0. \qqed
\end{align*}
\end{pf}
Now, we use Lemma~\ref{lem:small} to prove another combinatorial identity.
\begin{lem}\label{lem:big}
Let $k$ and $s$ be positive integers such that $r=k-s\geq 1$. Then
\begin{equation}\label{eqn:minus}
\sum_{i=0}^{k-s-1} (-1)^{i+1}{s+i \choose i}{k \choose s+i+1}=\sum_{i=0}^{r-1} (-1)^{i+1}{s+i \choose i}{s+r \choose s+i+1}=-1.
\end{equation}
\end{lem}
\begin{pf}
We use induction on $r=k-s$. Clearly, the result is true for $r=1$.
Assume that equation~(\ref{eqn:minus}) holds for $r-1$. Now, we prove equation~(\ref{eqn:minus}) for $r$.
Then
\begin{align*}
 \sum_{i=0}^{r-1} (-1)^{i+1}{s+i \choose i}{s+r \choose s+i+1}=\sum_{i=0}^{r-1} (-1)^{i+1}{s+i \choose i}\left({s+r-1 \choose s+i}+{s+r-1 \choose s+i+1}\right)=&\\
 \sum_{i=0}^{r-1} (-1)^{i+1}{s+i \choose i}{s+r-1 \choose s+i}
 +\sum_{i=0}^{r-1-1} (-1)^{i+1}{s+i \choose i}{s+r-1 \choose s+i+1}=-1,
\end{align*}
where the last equality follows from the induction hypothesis and 
Lemma~\ref{lem:small}. \qed
\end{pf}
Now, we can prove the following theorem.
\begin{thm}\label{thm:Jcong}
Let $\D$ be an $\oand$ such that  $k\geq s+1$ and $\ell \in [k-s]$.
Then for $u \subseteq [k]$ and $|u|=s+\ell$, $$J^{\x}_u(i_1,\ldots,i_k)=
\mu_u(i_1,\ldots, i_k)n^s,$$ where 
\begin{equation}\label{eqn:formula}
\mu_u(i_1,\ldots, i_k)
\equiv(-1)^\ell \lambda \binom{s+\ell-1}{\ell-1}\equiv(-1)^{|u|-s}\lambda \binom{|u|-1}{|u|-s-1}\quad ({\rm mod}\ n).
\end{equation}
\end{thm}
\begin{pf}
We prove this result by induction on $k$. For $k=s+1$, we have $\ell=1$. Then 
by equation~(\ref{eqn:Ju}) and Lemma~\ref{lem:ii} ($ii$) we have
\begin{equation*}
J^{\x}_u (i_1,\ldots,i_{s+1})= n^s(nx(i_1,\ldots,i_{s+1})-\lambda)  \quad \mbox{for } |u|=s+1.
\end{equation*}
So,
$J^{\x}_u(i_1,\ldots,i_{s+1})=\mu_u(i_1,\ldots, i_{s+1})n^s$, where $|u|=s+1$, and
\begin{equation*}
\mu_u(i_1,\ldots, i_{s+1})\equiv (-1) \lambda \binom{s+1-1}{1-1}\equiv-\lambda \quad (\text{mod}\   n).
\end{equation*}
Now, assume that the result is true for $s+1 \leq k<k'$ or equivalently true for $s+1 \leq |u|<k'$ (by Lemma~\ref{lem:projD}), and prove it for $k'$ or equivalently for $|u|=k'$ (by Lemma~\ref{lem:projD}). 
Equations~(\ref{eqn:unnumbered}) and~(\ref{eqn:Jr})  imply
\begin{equation}\label{eqn:referee3}
	\sum_{v \subseteq u}J_v^{\x}(i_1,\ldots,i_{k'})=n^{|u|}
	\sum_{\{i_j\, \mid \, j \notin u\}}x(i_1,\ldots,i_{k'}).
\end{equation}
For $|u|=k'$, 
by equation~(\ref{eqn:referee3}), Lemma~\ref{lem:ii}, and the induction hypothesis we have
\begin{equation*}
J^{\x}_u (i_1,\ldots,i_{k'})=n^s\left[n^{k'-s}\sum_{\{i_j\, \mid \, j \notin u\}} x(i_1,\ldots,i_{k'})-\lambda -\sum_{s+1\leq|\gamma|< k'}\mu_{\gamma}(i_1,\ldots, i_{k'})\right]\+, 
\end{equation*}
where for  $|\gamma|\in \{s+1,\ldots,k'-1\}$
\begin{equation*}
\mu_\gamma(i_1,\ldots, i_{k'})\equiv (-1)^{|\gamma|-s} \lambda \binom{|\gamma|-1}{|\gamma|-s-1} \quad (\text{mod}\ n).
\end{equation*}
Then,
\begin{align*}
\mu_u(i_1,\ldots, i_{k'})& \equiv  n^{k'-s}\sum_{\{i_j\, \mid \, j \notin u\}} x(i_1,\ldots,i_{k'})-\lambda -\sum_{s+1\leq|\gamma| < k'}\mu_\gamma(i_1,\ldots, i_{k'})  
&&   (\mbox{mod }n)\\
 &\equiv -\lambda -\sum_{s+1\leq|\gamma|< k'}\mu_\gamma(i_1,\ldots, i_{k'}) && (\mbox{mod }n) \\ 
& \equiv -\lambda - \lambda\left[\sum_{\ell=1}^{k'-s-1}(-1)^{\ell}{k'  \choose  s+\ell}{s+\ell-1  \choose  \ell-1}\right] &&
(\text{mod}\ n)\\
& \equiv   -\lambda-\lambda\left[\sum_{\ell=0}^{k'-s-2}(-1)^{\ell+1}{k'  \choose  s+\ell+1}{s+\ell \choose  \ell}\right] && (\text{mod}\ n).
\end{align*}
Now, by Lemma~\ref{lem:big}
\begin{align*}
\mu_u(i_1,\ldots, i_{k'})&\equiv -\lambda-\lambda\left[ -1+(-1)^{k'-s+1}{k'-1
 \choose  k'-s-1}{k' \choose  k'}\right] &&
(\text{mod}\ n)\\
&\equiv  \lambda (-1)^{k'-s}{k'-1 \choose  k'-s-1} &&
(\text{mod}\ n). \qqed
\end{align*}
\end{pf} 
\section{Decomposing $\mathbb{Q}^{\X}$ into irreducible subrepresentations under $G^{{\rm iso}}(k,n)$ and $G(k)^{\rm OD}$
}\label{sec:dim2}
In this section we determine the decomposition of $\mathbb{Q}^{\X}$ into irreducible $\mathbb{Q}$-subrepresentations under the action of the largest known subgroup of $G^{\text{LD(\ref{ilp:BF}})}$. 
 
 The following theorem follows easily from~\cite[p.~134]{Diaconis1988}.
 \begin{thm}\label{thm:double}
 Let $G$ be a finite group acting on a set $X$. Let $G$ act on $X^k$ by
 $$ 
g(x_1,\ldots, x_k)= (g x_1,\ldots, g x_k)$$ for $g \in G$,
 $(x_1,\ldots,x_k)\in X^k$, and  
let 
$$F(h)=|\{x\in X\mid hx=x\}|$$ for each $h \in G.$ Then the following hold.
\begin{enumerate}
\item For each $k \in \mathbb{Z}^{\geq1}$ $$\frac{1}{|G|}\sum_{h \in G} F(h)^k=|\{\text{orbits of $G$ on $X^k$}\}|.$$
\item Let $R : G \rightarrow {\rm Aut}_{\mathbb C}(\mathbb{C}^X)$ be the permutation $\mathbb{C}$-representation associated to $X$, i.e., for the standard basis  $\{e_x \mid x\in X\}$ of $\mathbb{C}^X$, $R(h)e_x=e_{hx}.$ Let
$$ \mathbb{C}^X= m_0V_0 \obot \cdots \obot m_bV_b $$ be the decomposition of $\mathbb{C}^X$ into irreducible $\mathbb{C}$-subrepresentations, where $m_i\geq 1$ is the  multiplicity of  $V_i$, i.e., $m_i$ is the number of times the irreducible $\mathbb{C}$-representation $V_i$ appears up to equivalence  in the decomposition $ \mathbb{C}^X= m_0V_0 \obot \cdots \obot m_bV_b $. 
Then 
$$ \sum_{i=0}^bm_i^2=|\{\text{orbits of $G$ on $X^2$}\}|.$$
\item The multiplicities $m_i$ satisfy   $m_i=1$ for each $i$ 
if and only if $b+1=|\{\text{orbits of $G$ on $X^2$}\}|$.
\end{enumerate}
 \end{thm}
 Let $\X$ be an $n^k\times k$ array, and $\{\text{rows of }\X\}$  consists of the  $n^k$ symbol combinations from $\{l_1,\ldots, l_n\}^k$ as in Section~\ref{sec:ANOVA}. Then $G^{\rm iso}(k,n)$ acts on the elements of 
$\{\text{rows of } \X\}$ via its action on the elements of $\{\text{columns of } \X\}$. Consequently $G^{\rm iso}(k,n)$ acts on the elements of
$\{e_{\x_1},\ldots, e_{\x_{n^k}}\}$, where $\x_i$ is the $i$th row of $\X$.
 Hence, the resulting action of 
$G^{\rm iso}(k,n)$ on $\mathbb{Q}^{\X}$ is a 
permutation $\mathbb{Q}$-representation.

Recall that for two rows $\x_1$ and $\x_2$ of $\X$,   
 the number of non-zero entries in $\x_1-\x_2$ is $d(\x_1,\x_2)$,  i.e., the {\em Hamming distance}  between  $\x_1$ and $\x_2$. 
We need the following two lemmas to find the decomposition of 
$\mathbb{Q}^{\X}$ into irreducible $\mathbb{Q}$-subrepresentations 
under the action of $G^{\rm iso}(k,n)$.
\begin{lem}\label{lem:Giso}
Let $G^{\rm iso}(k,n)$ act on
$\{\text{rows of } \X\} \times \{\text{rows of } \X\}$ as in Theorem~\ref{thm:double}
 via its action on the elements of $\{\text{rows of } \X\}$.
For $i=0,1,\ldots,k$, let $O_i\subset \{\text{rows of } \X\} \times \{\text{rows of } \X\}$ be such that 
$(\x_1,\x_2)\in O_i$ if and only if $d(\x_1,\x_2)=i$.
 Then the orbits of $G^{\rm iso}(k,n)$ on $\{\text{rows of } \X\} \times \{\text{rows of } \X\}$ are 
 $O_0,O_1,\ldots,O_k$.
 \end{lem}
 \begin{pf}
 First, $G^{\rm iso}(k,n)=\left(\prod_{i=1}^kS_{\{l_1,\ldots, l_n\}_i}\right)
 \rtimes S_{\{1,\ldots, k\}}$, where 
 $S_{\{l_1,\ldots, l_n\}_i}$  permutes the symbols $\{l_1,\ldots,l_n\}$ 
 in the $i$th column of $\X$, and $S_{\{1,\ldots, k\}}$ 
 permutes the columns of $\X$.
 Clearly, $d(\x_1,\x_2)=d(g\x_1,g\x_2)$ $\forall g \in G^{\rm iso}(k,n)$.
 Hence, $G^{\rm iso}(k,n)$ acts on the elements of each $O_i$. To show that $G^{\rm iso}(k,n)$ acts transitively on the elements of each $O_i$, let $(\x_1,\x_2)\in O_i$. 
 Since $G^{\rm iso}(k,n)$ acts transitively on the elements of $\{\text{rows of } \X\}$,
 there exists some $g_1 \in G^{\rm iso}(k,n) \, \ni \, g_1\x_1=\x'_1=(l_1,\ldots,l_1)$ and $g_1\x_2=\x'_2$, where 
 $d(\x'_1,\x'_2)=i$. 
 Then, there exists $g_2 \in G^{\rm iso}(k,n) \, \ni \, 
g_2\x'_2=(l_n,\ldots,l_n,l_1,\ldots,l_1)$. Hence, 
 $(g_2g_1\x_1,g_2g_1\x_2)=((l_1,\ldots,l_1),(l_n,\ldots,l_n,l_1,\ldots,l_1))$ for arbitrary $(\x_1,\x_2)\in O_i$. This proves that $G^{\rm iso}(k,n)$ acts transitively on the elements of each $O_i$.
 \qed
 \end{pf}
 \begin{lem}\label{lem:Gkiso}
 Let  $\{\text{rows of }\X\}$ consist of all $n^k$ combinations from
  $\{l_1,\ldots,l_n\}^k$, and $L_u$,  $L'_u$ be as in Theorem~\ref{thm:generalk}.   Then for each $r$ such that $0 \leq r\leq k$ each of the subspaces over $\mathbb{Q}$,  $\mathbb{C}$  $$U_r=\bigobot_{ u \subseteq [k],|u|=r} L_u, \quad U'_r=\bigobot_{ u \subseteq [k],|u|=r} L'_u $$ 
is invariant 
under the action of $G^{\rm iso}(k,n)$. Moreover, 
\begin{equation}\label{eqn:CX3}
	\mathbb{Q}^{\X}=\bigobot_{r=0}^k U_r,
\end{equation}
and
\begin{equation}\label{eqn:CX2}
	\mathbb{C}^{\X}=\bigobot_{r=0}^k U'_r.
\end{equation} 
 \end{lem}
 \begin{pf}
Let $U_{\varepsilon_j,j}$ for $j \in [k]$ be as in Theorem~\ref{thm:generalk}. Let $\vv \in \mathbb{Q}^{\X}$  be of the form $\vv=\vv_{\varepsilon_1,1} \otimes \cdots \otimes \vv_{\varepsilon_k,k}$, where $\vv_{\varepsilon_j,j} \in U_{\varepsilon_j,j}$. Then an element $\left((h_1, \ldots, h_k ), g\right) \in G^{\rm iso} (k, n)$ acts on $\vv$ by  
$$\left((h_1, \ldots, h_k ), g\right) (\vv_{\varepsilon_1,1} \otimes \cdots \otimes \vv_{\varepsilon_k,k}) 
= h_1 \vv_{\varepsilon_{g^{-1}(1)},g^{-1} (1)} \otimes \cdots \otimes h_k 
\vv_{\varepsilon_{g^{-1}(k)},g^{-1} (k)}.$$ 
Consequently, 
$$\left((h_1, \ldots, h_k ), g\right) (L_u)  = L_{g^{-1}u},$$
and it is immediate that $U_r$ is invariant under $G^{{\rm iso}}(k, n)$. 
Hence, each of the $k+1$ subspaces   $$U_r=\bigobot_{ u \subseteq [k],|u|=r} L_u$$
over $\mathbb{Q}$ for $r=0,1,\ldots,k$ with dim($U_r)={k \choose r}(n-1)^r$ 
is invariant under the action of $G^{\rm iso}(k,n)$. Equations~(\ref{eqn:CX3}) and~(\ref{eqn:CX2}) follow from equations~(\ref{eqn:all}) and~(\ref{eqn:all2}).

The proof for the $\mathbb{C}$ vector spaces $U'_r$ is obtained from the above proof by replacing $\mathbb{Q}$ with $\mathbb{C}$, $U_{\varepsilon_j,j}$ with $U'_{\varepsilon_j,j}$, and $L_u$ with $L'_u$.
\qed
 \end{pf}
 \begin{thm}\label{thm:Gkiso}
Decomposition~(\ref{eqn:CX2})~((\ref{eqn:CX3}))
   is the orthogonal decomposition of 
$\mathbb{C}^{\X}$  ($\mathbb{Q}^{\X}$) into $k+1$ 
irreducible $\mathbb{C}$-representations ($\mathbb{Q}$-representations)
under the action of $G^{\rm iso}(k,n)$.
 \end{thm}
 \begin{pf}
  Let
$ \mathbb{C}^{\X}= m_0V_0 \obot \cdots \obot m_bV_b $ be the decomposition of $\mathbb{C}^{\X}$ into irreducible subrepresentations, where $m_i\geq 1$ is the multiplicity of the irreducible $\mathbb{C}$-representation $V_i$. 
By Theorem~\ref{thm:double} and Lemma~\ref{lem:Giso},
 $ \sum_{i=0}^b m_i^2=k +1$.
 Since by Lemma~\ref{lem:Gkiso} $$k+1 \leq \sum_{i=0}^b m_i\leq \sum_{i=0}^b m^2_i=k+1,$$ we get  
 \begin{equation}\label{eqn:k+1}
 \sum_{i=0}^b m_i=\sum_{i=0}^b m^2_i=k+1.
 \end{equation}
 Then $\sum_{i=0}^b m_i=\sum_{i=0}^b m^2_i$ implies $m_i=1$ for $i=0,1,\ldots,k$. Hence, by equation~(\ref{eqn:k+1}) $k+1=b+1$, and for 
 $r=0,1,\ldots,k$  each  of the subspaces $U'_r$ in Lemma~\ref{lem:Gkiso}  is an irreducible $\mathbb{C}$-representation. 
 Now, the result for $\mathbb{Q}^{\X}$ follows from the
  result for $\mathbb{C}^{\X}$, 
  Lemma~\ref{lem:Gkiso}, and the fact that $\mathbb{Q}\subset \mathbb{C}$.
 \qed
 \end{pf}
   Let  $\{\text{rows of }\ZZ\}$  consist of all $2^k$ combinations from
  $\{-1,1\}^k$.
For $j \in [k]$ define the column operation $R_j$ on $\ZZ$ to be
{\small
\begin{equation}\label{eqn:MZZ}
\ZZ=\begin{bmatrix}
\z_1& \cdots &  \z_j&\cdots&\z_k\\
 \end{bmatrix} \quad  \stackrel{R_j}{\longrightarrow}\quad
 \begin{bmatrix}
\z_1\odot\z_j& \cdots&\z_{j-1}\odot\z_j& \z_j&\z_{j+1}\odot\z_j&\cdots & \z_k\odot \z_j\\
 \end{bmatrix}.
\end{equation}} 
Let 
\begin{equation}\label{eqn:Gkod}
G(k)^{\rm OD}=\langle R_1,\ldots,R_k, G^{\rm iso}(k,2)\rangle.
\end{equation} 
Then both $G^{\rm iso}(k,2)$ and $G(k)^{\rm OD}$ act on the rows of $\ZZ$. In this case,
\begin{equation}\label{eqn:Gkoddd}
G^{\rm iso}(k,2)=\left(\prod_{i=1}^kS_{\{-1,1\}_i}\right)
 \rtimes S_{\{1,\ldots, k\}},
\end{equation} 
 where $S_{\{-1, 1\}_i}$  swaps the symbols $\{-1, 1\}$ 
 in the $i$th column of $\ZZ$, and $S_{\{1,\ldots, k\}}$ 
 permutes the columns of $\ZZ$. Moreover, $R := \langle R_1, \ldots, R_k \rangle \cong
S_{k+1}$,  $R \cap G^{{\rm iso}} (k, 2) = S_{\{1,\ldots,k\}}$, and $\prod_{i=1}^kS_{\{-1,1\}_i}$ is still 
normal in $G(k)^{{\rm OD}}$~\cite{Arquette2016}. Thus $G(k)^{{\rm OD}} = R\left(\prod_{i=1}^kS_{\{-1,1\}_i}\right) = \left(\prod_{i=1}^kS_{\{-1,1\}_i}\right)R \cong (S_2^k ) \rtimes S_{k+1}$. 
 \begin{lem}\label{lem:contains2}
Let $n=2$ in ILD~(\ref{ilp:BF}), and 
 $\x \in \mathbb{Z}^{2^k}$ be such that $x(i_1, \ldots, i_k)$  is the number of times the  symbol combination $(i_1,\ldots,i_k)$ with $(i_1,\ldots,i_k)^{\top} \in \{-1,1\}^k$ appears as a row of a sought after $\oant$ with symbols from $\{-1,1\}$. 
 Let $G(k,2,s)^{\rm LD}$ be the symmetry group of the LD relaxation of  ILD~(\ref{ilp:BF}). Then, $G(k,2,s)^{\rm LD} \geq G(k)^{\rm OD}$  
if and only if $s$ is even. Hence, for even $s$, $|G(k,2,s)^{\rm LD}|\geq |G(k)^{\rm OD}|=|S_2^k \rtimes S_{k+1}|=(k+1)!2^k$.
\end{lem}
\begin{pf}
The proof follows from the proof of Lemma 11 in~\cite{Geyer2018}. \qed
\end{pf}
Next, we determine the orbits of  $G^{\rm iso}(k,2)$ and $G(k)^{\rm OD}$ 
 on $\{\text{rows of }\ZZ\}\times\{\text{rows of } \ZZ\}$.
 \begin{lem}\label{lem:OOO}
  Let $O_0=O'_0\subset \{\text{rows of }\ZZ\}\times\{\text{rows of } \ZZ\}$ be such that $O_0=O'_0=\cup_{\z \in \{\text{rows of }\ZZ\}}\{(\z,\z)\}$, and
for $i \in [k]$, let $O'_i\subset \{\text{rows of }\ZZ\}\times\{\text{rows of } \ZZ\}$ be such that 
$(\z_1,\z_2)\in O'_i$ if and only if $d(\z_1,\z_2)=i$ or 
$d(\z_1,\z_2)=k+1-i$.
 Then the orbits of $G(k)^{\rm OD}$ on $\{\text{rows of }\ZZ\}\times\{\text{rows of } \ZZ\}$ are 
 $O_0=O'_0,O'_1,\ldots,O'_{\lceil k/2\rceil}$.
 \end{lem}
 \begin{pf}
 Clearly, $O_0=O'_0=\cup_{\z \in \{\text{rows of }\ZZ\}}\{(\z,\z)\}$ is an orbit of $G(k)^{\rm OD}$.
 Let  $$O_i = \{(\z_1, \z_2 ) \mid d(\z_1, \z_2) = i\}.$$ Then by
Lemma~\ref{lem:Giso},  $O_i$ is an orbit of $G^{\rm iso}(k,2)=\left(\prod_{i=1}^kS_{\{-1,1\}_i}\right) 
 \rtimes S_{\{1,\ldots, k\}}$ on $\{\text{rows of } \ZZ\} \times \{\text{rows of } \ZZ\}$. 
 By the
definition of $O_i'$, it is trivial that 
$O_i' = O_i \cup O_{k+1-i}$. Let $R_j$ for $j \in [k]$ be as in equation~(\ref{eqn:MZZ}). As $O_i$ is an orbit of 
$G^{{\rm iso}} (k, 2) \leq G(k)^{\rm{OD}}$, and since
$d(R_j \z_1, R_j \z_2 ) = k + 1 - d(\z_1, \z_2)$ if and only if  ${z_1}_j \neq {z_2}_j$, 
it follows that $O_i \cup O_{k+1-i}$ is a
$G(k)^{\rm{OD}}$-orbit for $i=1,\ldots,\lceil k/2\rceil$.
 \qed
 \end{pf}
 The proof of the following lemma mimics the proof of Lemma 7 
 in~\cite{Geyer2018}.
 \begin{lem}\label{lem:Jr}
  Let the rows of $\ZZ$ be all $2^k$  
  vectors in $\{-1,1\}^k$, and  $\x \in \mathbb{C}^{\ZZ}$.
Let   $u \subseteq [k]$ be such that $|u|=r\geq0$, and $G(k)^{\rm OD}$ be as in equation~(\ref{eqn:Gkod}).
Let $g \in G(k)^{\rm OD}$ and $g(\x)$ 
be obtained after $g$ is applied to $\x$.
 Then $$J_r(u)^{g(\x)}=\pm J_{r'}(u')^{\x}$$ for some $u' \subseteq [k]$, where 
 \begin{eqnarray}\label{eqn:spann}
|u'|=r'=\left\{
\begin{array}{rl}
\text{$r$ or $r+1$} & \quad \text{if $r>0$ and $r$ is odd,}\\
\text{$r$ or $r-1$} &\quad  \text{if $r>0$ and $r$ is even,}\\
\text{$0$} &\quad  \text{if $r=0$.}
\end{array}\right.
\end{eqnarray}
 \end{lem}
 \begin{pf}
Since each $g \in G(k)^{\rm OD}$ permutes the rows of $\ZZ$, 
$G(k)^{\rm OD}$ acts on  $\mathbb{C}^{\ZZ}$ and the resulting representation of $G(k)^{\rm OD}$ is a permutation $\mathbb{C}$-representation.
For each $i\in[k]$ let $R_i$ be defined as in equation~(\ref{eqn:MZZ}). 
Then, 
\begin{align}\label{eqn:Ri}
J_r(u)^{R_i(\x)}= \begin{cases}
J_r(u)^{\x} & \text{if $r>0$, $r$ is even, and $i \notin u$,}\\
J_{r-1}(u\backslash\{i\})^{\x} & \text{if $r>0$, $r$ is even, and $i \in u$,}\\
J_{r+1}(u\cup \{i\})^{\x} & \text{if $r>0$, $r$ is odd, and $i \notin u$,}\\
J_r(u)^{\x} &\text{if $r>0$, $r$ is odd, and $i \in u$,}\\
J_r(u)^{\x}=\sum_{i_1,\ldots,i_k}x(i_1,\ldots,i_k) &\text{if $r=0$.}
\end{cases}
\end{align}
Let $R=\langle R_1,\ldots,R_k \rangle$, and $\prod_{i=1}^kS_{\{-1,1\}_i}$ be the group of all possible sign switches of columns of $\ZZ$.
 Then by the proof of Lemma 4 in~\cite{Geyer2018}, $g=g_1g_2$
 for each $g \in G(k)^{\rm OD}$, where $g_1 \in R$   and $g_2 \in \prod_{i=1}^kS_{\{-1,1\}_i}$. 
Hence, by equation~(\ref{eqn:Ri}),  
\begin{equation*}\label{eqn:isomorphism}
J_r(u)^{g(\x)}=J_r(u)^{g_1(g_2(\x))}=J_{r'}(u')^{g_2(\x)}
\end{equation*}
for some $u' \subseteq [k]$, and $r'=|u'|$ is as in equation~(\ref{eqn:spann}). 
 Now, $g_2(\x)$ is obtained by permuting the rows of $\ZZ$ that corresponds to multiplying a subset of columns of $\ZZ$ by $-1$. Therefore,
 \begin{equation*}
 J_r(u)^{g(\x)}=J_{r'}(u')^{g_2(\x)}=\pm J_{r'}(u')^{\x}.\qqed
 \end{equation*}
 \end{pf}
The next two lemmas will be used to find the decomposition of $\mathbb{Q}^{\X}$ into irreducible $\mathbb{Q}$-subrepresentations under the action of $G(k)^{\rm OD}$.
  \begin{lem}\label{lem:Jrirred}
Let $G(k)^{\rm OD}$ be as in equation~(\ref{eqn:Gkod}), 
$U'_r$ be as in Lemma~\ref{lem:Gkiso}.
Let $W'_0=U'_0$, $W'_j=U'_{2j-1}\obot U'_{2j}$
for $j \in [\lceil k/2\rceil-1]$, and
$$
W'_{\lceil \frac{k}{2}\rceil}=\begin{cases}
U'_{k-1}\obot U'_{k} \quad \, \text{if } k  \text{ is even,} \\
U'_{k}\quad \quad \quad \quad \, \, \, \text{otherwise.}\\
\end{cases}$$
Then, 
\begin{equation}\label{eqn:allall}
 ({\mathbb C^2})^{\otimes k}\cong \mathbb{C}^{\ZZ} = 
  \bigobot_{j=0}^{\lceil \frac{k}{2}\rceil}W'_j
\end{equation} 
is an orthogonal decomposition of $\mathbb{C}^{\ZZ}$ into  invariant subspaces under the action of $G(k)^{\rm OD}$.
 \end{lem}
 \begin{pf}
Let $\HH$  be as in equation~(\ref{eqn:M}). Then by the invertibility of $(\HH^{\top})/2^k$ 
$$\mathbb{C}^{\ZZ}=\text{Col}_{\mathbb{C}}(\frac{\HH^{\top}}{2^k}).$$ 
 Let 
 $$\h_{\{i_1,\ldots,i_{j}\}}=
 \frac{1}{2^k}\z_{i_1,\ldots,i_j},$$ where 
  and $\z_{i_1,\ldots,i_j}$ is as in equation~(\ref{eqn:M}). 
 Then $$\bigcup_{j=0}^k\bigcup_{i_1<\cdots <i_{j}}\{\h_{\{i_1,\ldots,i_{j}\}}\}=\{{\rm Columns \ of} \left(\frac{\HH^{\top}}{2^k}\right)\}.$$
Let  $\mathcal{B}_0=\{\h_{\emptyset}$\}, and for $j \in [\lceil k/2\rceil-1]$
let $$\mathcal{B}_j=\bigcup_{i_1<\cdots <i_{2j}}\{\pm \h_{\{i_1,\ldots,i_{2j}\}}\}\bigcup \bigcup_{i_1<\cdots <i_{2j-1}}\{\pm
 \h_{\{i_1,\ldots,i_{2j-1}\}}\},$$ and 
$$
\mathcal{B}_{\lceil \frac{k}{2}\rceil}=\begin{cases}
\bigcup_{i=1}^k \{\pm\h_{[k]\backslash \{i\}}\}\bigcup
\{\pm\h_{[k]}\} \quad \text{if } k \text{ is even,} \\
\{\pm\h_{[k]}\}    \hspace{3.9cm} \text{otherwise.}\\
\end{cases}
$$ Now, equations~(\ref{eqn:MZZ}), (\ref{eqn:Gkod}), and (\ref{eqn:Gkoddd})  imply that each element of  $G(k)^{\rm OD}$ acts on the elements of $\mathcal{B}_j$ for $j=0,1,\ldots,$ $\lceil k/2\rceil$ 
 as a signed permutation (a permutation that may or may not be followed by sign changes). So, Span$_{\mathbb{C}}(\mathcal{B}_j)=W'_j$ is invariant under the action of  $G(k)^{\rm OD}$.
 \qed
 \end{pf} 
 The following lemma follows from the fact that $\h_{\{i_1,\ldots,i_{j}\}}\in \mathbb{Q}^\ZZ$  
 $\, \forall  \{i_1,\ldots,$ $i_{j}\} \, \ni \, i_1<\cdots<i_j$.
 \begin{lem}\label{lem:valid}
Lemma~\ref{lem:Jrirred} remains valid if the field of scalars $\mathbb{C}$
 is replaced with $\mathbb{Q}$.
 \end{lem} 
  \begin{thm}\label{thm:GkOD}
 Let $G(k)^{\rm OD}$ be as in Lemma~\ref{lem:Jrirred}. Then  decomposition~(\ref{eqn:allall}) in Lemma~\ref{lem:Jrirred} is the orthogonal decomposition of $\mathbb{C}^{\ZZ}$  into irreducible $\mathbb{C}$-subrepresentations.
 \end{thm} 
 \begin{pf}
 Let $$ \mathbb{C}^{\ZZ}= m_0V_0 \obot \cdots \obot m_bV_b $$
 be the decomposition of $ \mathbb{C}^{\ZZ}$ into 
 irreducible $\mathbb{C}$-representations under the action of $G(k)^{\rm OD}$, where $m_i\geq 1$ is the multiplicity of the $\mathbb{C}$-representation $V_i$.
 By Lemma~\ref{lem:Jrirred}, $\lceil k/2\rceil+1\leq \sum_{i=0}^b m_i$.
 Moreover, by Theorem~\ref{thm:double} and Lemma~\ref{lem:OOO} 
 $$
 \left \lceil \frac{k}{2}\right \rceil+1=\sum_{i=0}^b m_i^2\geq\sum_{i=0}^b m_i  \geq \left \lceil \frac{k}{2}\right \rceil+1.
 $$
 Hence, $m_i=1$ for $i=0,1,\ldots,\lceil k/2\rceil+1$, and $\lceil k/2\rceil+1= b+1$. \qed
  \end{pf} 
  \begin{cor}\label{cor:GkODQ}
 Let $W_j$ be obtained in the same way $W'_j$ is obtained in Lemma~\ref{lem:Jrirred} after replacing $U'_{r}$ with $U_{r}$ in Lemma~\ref{lem:Gkiso}. Then
 $$ ({\mathbb Q^2})^{\otimes k}\cong \mathbb{Q}^{\ZZ} = 
  \bigobot_{j=0}^{\lceil \frac{k}{2}\rceil}W_j $$
is the orthogonal decomposition of $\mathbb{Q}^{\ZZ}$  into irreducible $\mathbb{Q}$-subrepresentations.
  \end{cor}
  \begin{pf}
  Invariance of the subspaces $W_j$ under the action of $G(k)^{\rm OD}$ 
  follows from  Lemma~\ref{lem:valid}.
  Irreducibility of $W_j$ follows from Theorem~\ref{thm:GkOD} and the fact that $\mathbb{Q}\subset \mathbb{C}$.\qed
\end{pf} 
\section{Decreasing the number of possible values of dim($P^{(k,s,\lambda)}_{n;I}$)}\label{sec:dim}
In this section, by using representation theory, we drastically decrease the number of possible values of dim($P^{(k,s,\lambda)}_{n;I}$). We also determine the corresponding sets of potentially valid equality constraints for Aff($P^{(k,s,\lambda)}_{n;I}$). These  are the only sets of equality constraints up to equivalence that can be implied by the integrality constraints of ILD~(\ref{ilp:BF}). 
First, we show that the ILD 
\begin{eqnarray}\label{ilp:J}
\begin{array}{rl} 
 & J^{\x}_{\emptyset}(i_1,\ldots, i_k)=\lambda n^s,\  \\
  & J^{\x}_{u}(i_1,\ldots, i_k)=0,\quad
   \forall  (i_1,\ldots,i_k)^{\top} \in \{l_1,\ldots,l_n\}^k, \ \ \text{for } u \ni 1\leq |u|\leq s,\\
 &    0 \leq x(i_1,\ldots, i_k) \leq  p_{\max}, \quad x(i_1,\ldots, i_k) \in 
 \mathbb{Z}, \quad \forall   (i_1,\ldots, i_k)^{\top} \in 
 \{l_1,\ldots, l_n\}^k  
 \end{array}
\end{eqnarray}
 is equivalent to ILD~(\ref{ilp:BF}), where $p_{\max}$ is computed as in ILD~(\ref{ilp:BF}).
\begin{lem}\label{lem:BFJ}
 The equality constraints of LD~(\ref{ilp:BF}) can be obtained as linear combinations of the  equality constraints of LD~(\ref{ilp:J}) and vice versa.
\end{lem} 
\begin{pf}
First, by equations~(\ref{eqn:grandmean}) and (\ref{eqn:Jr}) 
\begin{equation}\label{eqn:xempty}
J^{\x}_{\emptyset}(i_1,\ldots, i_k) =n^kx_{\emptyset} (i_1, \ldots, i_k)= 
\sum_{i_1, \ldots, i_k} x(i_1, \ldots, i_k)=\lambda n^s,
\end{equation} 
and
\begin{equation}\label{eqn:xu}
J^{\x}_{u}(i_1,\ldots, i_k)=n^kx_u (i_1, \ldots, i_k) = 
n^{|u|}
\sum_{ \{i_j \,\mid \, j \not\in u\}}  
x(i_1, \ldots, i_k) 
- n^k \sum_{ v \subsetneq u } x_v (i_1, \ldots, i_k). 
\end{equation}
Moreover, equality constraints of LD~(\ref{ilp:BF}) imply
\begin{equation}\label{eqn:sumu}
\sum_{ \{i_j \,\mid \, j \not\in u\}}  
x(i_1, \ldots, i_k)=\lambda n^{s-|u|} \quad \text{for $u \ni |u| \leq s$}. 
\end{equation}
Combining equations~(\ref{eqn:xempty}),~(\ref{eqn:xu}), and~(\ref{eqn:sumu}) we get
$$ J^{\x}_{u}(i_1,\ldots, i_k)=0\quad \text{ for each $ u \ni 1 \leq |u| \leq s$}.$$
Conversely, let  the  equality constraints of LD~(\ref{ilp:J}) hold. Then
these constraints and equation~(\ref{eqn:xu}) imply equations~(\ref{eqn:sumu}).
We conclude the proof by observing that equations~(\ref{eqn:sumu}) for $|u|=s$ are the equality constraints of LD~(\ref{ilp:BF}). \qed
\end{pf}
Both  ILD~(\ref{ilp:BF}) and ILD~(\ref{ilp:J}) have the same inequality constraints. Hence, by Lemma~\ref{lem:BFJ}  the LD relaxation feasible sets of ILD~(\ref{ilp:BF}) and ILD~(\ref{ilp:J}) are  the same. Consequently, the 
  feasible sets of ILD~(\ref{ilp:BF}) and ILD~(\ref{ilp:J}) are the same, and consist of the frequency vectors of all $\oand$. 
 The LD relaxation of ILD~(\ref{ilp:BF}) has $\sum_{j=0}^{s}{k\choose j}(n-1)^j$ non-redundant equality constraints~\cite{Rosenberg1995}. So, the dimensions of the feasible sets of both LDs~(\ref{ilp:BF}) and~(\ref{ilp:J}) are $n^k-\sum_{j=0}^{s}{k\choose j}(n-1)^j$. 

  For 
   an arbitrary function $y(i_1,\ldots,i_k) \in \mathbb{C}^{\X}$, define $J^{\y}_{u}(i_1,\ldots, i_k)$ by replacing $x(i_1,\ldots,i_k)$ with $y(i_1,\ldots,i_k)$ in 
  equations~(\ref{eqn:ANOVA}-\ref{eqn:Jr}). Let $S_n \wr S_k$  act on 
feasible points as described in equations~(\ref{eqn:act}).
Now, the following lemma is used to show that the action of $S_n \wr S_k$  drastically decreases the number of   possible values of $\mbox{dim(Conv}(S_n \wr S_k \,\,\x))$.
\begin{lem} \label{lem:u''}
If for each feasible point $\x$ of  ILD~(\ref{ilp:J}) 
 $$ 
 J^{\x}_{u'}(i_1,\ldots, i_k)=0 \quad\text{for some $u'\subseteq[k]\, \ni\, $  } |u'|\geq s+1, $$  $\forall(i_1,\ldots, i_k)^{\top} \in 
 \{l_1,\ldots, l_n\}^k$, then for $\y=\x- \lambda n^s/n^k\1_{n^k}$ we must also have 
 $$ J^{\y}_{u''}(i_1,\ldots, i_k)=J^{\x}_{u''}(i_1,\ldots, i_k)=0 \quad\forall u''\subseteq[k] \, \ni \, |u''|=|u'|, $$ $\forall (i_1,\ldots, i_k)^{\top} \in 
 \{l_1,\ldots, l_n\}^k$,   and all feasible points $\x$ of ILD~(\ref{ilp:J}).
 \end{lem}
 \begin{pf}
 First, $ J^{\y}_{u}(i_1,\ldots, i_k)=J^{\x}_{u}(i_1,\ldots, i_k)$ $\forall u\subseteq [k]\, \ni \, u\neq \emptyset$ since $\y=\x+a\1_{n^k}$ for some constant $a \in \mathbb{C}$. Now, since, $$J^{((h_1,\ldots,h_k),g)(\x)}_{u'}(i_1,\ldots, i_k)=J^{\x}_{g^{-1}(u')}(h_1(i_1),\ldots, h_k(i_k))$$ we get that
 $\left(\prod_{j=1}^kS_{\{\l_1,\ldots, l_n\}_j}\right)\rtimes S_{\{1,\ldots,k\}}$ acts transitively on the elements of $$\{J^{\x}_{u'}(i_1,\ldots, i_k)\ |\ |u'|=r\}$$ while preserving the feasible points of ILD~(39). Hence, 
 $$ J^{\y}_{u''}(i_1,\ldots, i_k)=J^{\x}_{u''}(i_1,\ldots, i_k)=0 \quad\forall u''\subseteq[k] \, \ni \, |u''|=|u'|, $$ $\forall (i_1,\ldots, i_k)^{\top} \in 
 \{l_1,\ldots, l_n\}^k$, and all feasible points $\x$ of ILD~(\ref{ilp:J}). 
  \qed 
\end{pf}
The following lemma strengthens  Lemma~\ref{lem:u''} when $n=2$,  
 $s$ is even, and $|u'|$ is even.
\begin{lem} \label{lem:u''2}
Let $n=2$, and $s$ be even. Let  $u'\subseteq[k]$ be such that 
 $|u'|\geq s+1$, and $|u'|$ be even. 
If for each feasible point $\x$ of  ILD~(\ref{ilp:J}) 
 $$ J^{\x}_{u'}(i_1,\ldots, i_k)=0  \quad 
\forall (i_1,\ldots, i_k)^{\top} \in \{l_1,\ldots,l_n\}^k, $$ then we must also have 
 $$ J^{\y}_{u''}(i_1,\ldots, i_k)=J^{\x}_{u''}(i_1,\ldots, i_k)=0 \quad \forall  u''\subseteq[k] \, \ni \, |u'|-1 \leq |u''|\leq |u'|, $$ $\forall (i_1,\ldots, i_k)^{\top} \in \{l_1,\ldots,l_n\}^k,$ all feasible points $\x$ of ILD~(\ref{ilp:J}), and  $\y=\x- \lambda n^s/n^k\1_{n^k}$.
 \end{lem}
 \begin{pf}
 First, $ J^{\y}_{u}(i_1,\ldots, i_k)=J^{\x}_{u}(i_1,\ldots, i_k)$ $\forall u\subseteq [k]\, \ni \, u\neq \emptyset$. Then the result follows, because  by Lemmas~\ref{lem:consistent} and~\ref{lem:Jr} for each $r \in [\lfloor k/2\rfloor]$, $G(k)^{\text{OD}}$ acts transitively on the elements of $$\{J^{\x}_{u'}(i_1,\ldots, i_k)\ |\ |u'|=2r\text{ or } |u'|=2r-1\}$$
 while preserving the feasible points of ILD~(\ref{ilp:J}). 
  \qed 
\end{pf}
For a given feasible point $\x$ of ILD~(\ref{ilp:J}), we will prove a theorem providing a restriction for all possible values of $\mbox{dim(Conv}(S_n \wr S_k \, \,\x))$  as well as the corresponding sets of equality constraints. However, we first need the following well-known lemma and the subsequent lemma. 
\begin{lem}\label{lem:dimgengroup}
Let $S$ be a set of vectors in $\mathbb{F}^n$, where $\mathbb{Q} \subseteq \mathbb{F} \subseteq\mathbb{R}$. Then
${\rm Aff}_{\mathbb{F}}(S)={\rm Span}_{\mathbb{F}}
(S)$ if and only if $\0_n \in {\rm Aff}_{\mathbb{F}}(S)$.
\end{lem}
\begin{lem}\label{lem:dimdim}
Let $\x$ be a feasible point of ILD~(\ref{ilp:J}), 
and $\y=\x-\lambda n^{s-k} \1_{n^k}$. 
 Then $${\rm dim(Conv}(S_n \wr S_k \, \,\y))={\rm dim(Span}
(S_n \wr S_k\, \, \y)).$$ 
\end{lem}
\begin{pf}
Let $G=S_n \wr S_k$.
Pick a feasible $\x$, then $\x\neq \0_{n^k}$. Let
 $$\z=\frac{1}{|G|}\sum_{g \in G}g\y= \frac{1}{|G|}\sum_{g \in G}g\x-\lambda n^{s-k} \1_{n^k}.$$
Then $g\z=\z$ for each $g \in G$.
Moreover, since $G$ acts on the set of indices of $\z$ transitively and $g\z=\z$ for each $g \in G$, $\z \in {\rm Conv}(G\,\,\y) \subseteq {\rm Aff}(G\,\,\y)$  must have equal coordinates. 
Now,  by Lemma~\ref{lem:dimgengroup}  it suffices to show that $\z=\0_{n^k}$ to finish the proof as $\z \in {\rm Aff}(G\,\,\y)$ and
$${\rm dim(Conv}(G\,\, \y))={\rm dim(Aff}(G\,\,\y)).$$ 
Observe that $\z+\lambda n^{s-k} \1_{n^k}$ is an equal coordinate vector in ${\rm Conv}(G\,\,\x)$.
Hence $\z+\lambda n^{s-k} \1_{n^k}$ is an equal coordinate vector satisfying all the constraints of LD~(\ref{ilp:J}).
However, the only such equal coordinate vector is $\x=\lambda n^{s-k} \1_{n^k}$. Hence, $\z=\0_{n^k}$. \qed
\end{pf}
\begin{thm}\label{thm:dim}
 Let $\x$ be a feasible point of ILD~(\ref{ilp:J}), and 
$\Omega_2$ be as in Theorem~\ref{thm:dimfin}. 
  Then the following hold.
 \begin{enumerate} 
 \item[(i)] There exists some  $T\subseteq \Omega_2$ 
 such that $${\rm dim(Conv}\left(S_n \wr S_k\, \, \x)\right)=n^k-\sum_{j=0}^s{k \choose j}(n-1)^j-\sum_{d\in T }{k \choose d}(n-1)^{d}.$$  
\item[(ii)] The equality constraints of ILD~(\ref{ilp:J}) together 
with the distinct equalities in 
\begin{equation}\label{eqn:convexhullGiso}
 J^{\x}_{u''}(i_1,\ldots, i_k)=0 \quad\forall  u''\subseteq [k] \text{ with } |u''| \in T,\quad \forall \text{$(i_1,\ldots, i_k)\in \{l_1,\ldots,l_n\}^k$}
\end{equation} determine 
${\rm Aff}({\rm Conv}(S_n \wr S_k \, \,\x))=
{\rm Aff}(S_n \wr S_k \, \, \x)$. 
\end{enumerate}
\end{thm}
\begin{pf}
Let $\y=\x-\lambda n^{s-k}\1_{n^k}$. To prove part ($i$), it suffices to show that 
\begin{equation*}
{\rm dim}{\rm (Span}(S_n \wr S_k\,\,\y))=n^k-\sum_{j=0}^s{k \choose j}(n-1)^j-\sum_{d \in T }{k \choose d}(n-1)^{d}
\end{equation*}
 for some $T\subseteq \Omega_2$ as  
 $${\rm dim}{\rm (Conv}(S_n \wr S_k \,\,\y))=
{\rm dim}{\rm (Conv}(S_n \wr S_k \,\,\x)), \quad {\rm dim}{\rm (Conv}(S_n \wr S_k\,\, \y))={\rm dim}{\rm(Span}
(S_n \wr S_k \,\,\y)),$$ 
where the second equation follows from Lemma~\ref{lem:dimdim}. 

First,
let $\A'\x=\bb'$ be the equality constraints of ILD~(\ref{ilp:J}). 
Then $\y \in \mbox{Null}_{\mathbb{Q}}(\A')$ as $\lambda n^{s-k}\1_{n^k}$ is a solution of $\A'\x=\bb'$. Moreover, since $\pi \x$ is a solution of $\A'\x=\bb'$ for each $\pi \in S_n \wr S_k$, $\pi \y \in \mbox{Null}_{\mathbb{Q}}(\A')$  for each $\pi \in S_n \wr S_k$.
Hence, 
 ${\rm Span}_{\mathbb{Q}}(S_n \wr S_k \,\,\y) \subseteq \mbox{Null}_{\mathbb{Q}}(\A')$, and consequently we get 
$ 
	\mbox{dim}\mbox{(Span}(S_n \wr S_k\,\, \y))=\mbox{dim}_{\mathbb{Q}}\mbox{(Span}_{\mathbb{Q}}(S_n \wr S_k \,\,\y)) \leq n^k-\sum_{j=0}^s{k \choose j}(n-1)^j.
$ 
Now, since $\mbox{Span}_{\mathbb{Q}}(S_n \wr S_k\,\, \y)\subseteq ({\mathbb Q}^n)^{\otimes k}$ is invariant under the action of $S_n \wr S_k$,  
$\mbox{Span}_{\mathbb{Q}}(S_n \wr S_k \,\,\y)$ in $\mathbb{Q}^{\X}$ must be an orthogonal direct sum of the irreducible subspaces in  decomposition~(\ref{eqn:CX3}). Hence, 
 if \begin{equation}\label{eqn:ayni}
 \mbox{dim}\mbox{(Span}(S_n \wr S_k\,\, \y))=\mbox{dim}_{\mathbb{Q}}\mbox{(Span}_{\mathbb{Q}}(S_n \wr S_k \,\, \y)) < n^k-\sum_{j=0}^s{k \choose j}(n-1)^j,
 \end{equation}
then $\mbox{Span}_{\mathbb{Q}}(S_n \wr S_k\,\, \y)$ must be orthogonal to at least one of the irreducible invariant subspaces $U_t$
in  decomposition~(\ref{eqn:CX3}) for some $t\geq s+1$. This implies that there exists a maximal collection of distinct $u_1,\ldots, u_r \subseteq [k]$ such that $|u_j|=d_j \leq k$, and
\begin{equation}\label{eqn:mustbezero}
J^{\y}_{u_j}(i_1,\ldots, i_k) = 0\quad \forall (i_1,\ldots, i_k)^{\top}\in \{l_1,\ldots,l_n\}^k \text{ and } j \in [r].
\end{equation}
On the other hand, based on the definition of $J_u^{\y}(i_1,\ldots,i_k)$
as a function of $\y$ it is easy to see that
\begin{align}\label{ilp:Jyy}
   J^{\y}_{\emptyset}(i_1,\ldots, i_k)&=0,\nonumber\\
   J^{\y}_{u}(i_1,\ldots, i_k)&=J^{\x}_{u}(i_1,\ldots, i_k)\quad\forall u \subseteq [k] \, \ni \, u \neq \emptyset. 
\end{align}
Hence, by 
equations~(\ref{eqn:mustbezero}) and~(\ref{ilp:Jyy}),
 we also have
 \begin{equation}\label{eqn:Jxx}
J^{\x}_{u_j}(i_1,\ldots, i_k) = J^{\y}_{u_j}(i_1,\ldots, i_k) = 0\quad \forall (i_1,\ldots, i_k)^{\top}\in \{l_1,\ldots,l_n\}^k  \text{ and } j \in [r]. 
 \end{equation}
Now, by Theorem~\ref{thm:Jcong},  $$ J^{\y}_{u}(i_1,\ldots, i_k)=J^{\x}_{u}(i_1,\ldots, i_k)=\mu_u(i_1,\ldots, i_k)n^s\quad\forall u \subseteq [k]\, \ni \, |u|\geq s+1, $$ where
$$\mu_u(i_1,\ldots, i_k)\equiv(-1)^{|u|-s} \lambda \binom{|u|-1}{|u|-s-1} \quad (\text{mod}\ n),$$ 
and $u \subseteq [k]$ with $|u| \geq s+1$.
Hence, if $n \ndiv \lambda \binom{|u|-1}{|u|-s-1}$ for some
 $u\subseteq [k]\, \ni \, |u|\geq s+1$, then
$$\mu_u(i_1,\ldots, i_k)\neq 0, 
$$
and $$ J^{\y}_{u}(i_1,\ldots, i_k)\neq 0. 
$$ 
Thus, $u_1,\ldots,u_r$ in equation~(\ref{eqn:mustbezero}) must be such that $|u_j|=d_j\in \{s+1,\ldots,k\} $ and $$\lambda \binom{d_j-1}{d_j-s-1}\equiv 0 \quad \mbox{(mod } \ n)$$
 for $j \in [r]$. Let $T=\{d_{j_1},\ldots,d_{j_{r'}}\}$ be the set of all distinct $d_j$ for $j=1,\ldots,r.$
Now by Lemma~\ref{lem:u''} and equation~(\ref{eqn:Jxx}), 
  $$ J^{\y}_{u''}(i_1,\ldots, i_k)=J^{\x}_{u''}(i_1,\ldots, i_k)=0 \quad\forall  u''\subseteq[k] \ni |u''| \in T,\quad \forall (i_1,\ldots, i_k)^{\top}\in \{l_1,\ldots,l_n\}^k, $$and  all feasible points $\x$ of ILD~(\ref{ilp:J}). This establishes equations~(\ref{eqn:convexhullGiso}).
   Hence,  each distinct $d_{j_i}\in T$  reduces $$\mbox{dim}_{\mathbb{Q}}\mbox{(Span}_{\mathbb {Q}}(S_n \wr S_k \,\,\y))=\mbox{dim(Span}(S_n \wr S_k\,\, \y))=\mbox{dim(Conv}(S_n \wr S_k\,\, \y))$$ by 
$\text{dim}(U_{d_{j_i}})={k \choose d_{j_i}}(n-1)^{d_{j_i}}$. This proves~($i$). Now, the equality constraints of 
ILD~(\ref{ilp:J}) together 
with equations~(\ref{eqn:convexhullGiso}) determine   
$\text{Aff}({\rm Conv}(S_n \wr S_k\,\, \x))=\text{Aff}(S_n \wr S_k\,\, \x)$. This proves~($ii$). 
 \qed 
\end{pf}
The following  corollary  is an immediate  consequence of 
Theorem~\ref{thm:dim} ($i$). 
\begin{cor}\label{cor:ndiv}
Let $k>s$, $n \ndiv \lambda \binom{s+\ell-1}{\ell-1}$ for $\ell \in [k-s]$, and $P^{(k,s,\lambda)}_{n;I}\neq \emptyset$. Then $${\rm dim}(P^{(k,s,\lambda)}_{n;I})= {\rm dim}(P^{(k,s,\lambda)}_{n})=n^k-\sum_{j=0}^s{k \choose j}(n-1)^j.$$
\end{cor}
Corollary~\ref{cor:ndiv} implies all the values of $\text{dim}(P^{(k,s)}_{n;I})$ with $k>s$ and $s \neq 0$ in Table~\ref{tab:PnI}. For each of these cases 
$\mbox{dim}(P^{(k,s)}_{n;I})=\mbox{dim}(P^{(k,s)}_{n})$ whenever $P^{(k,s)}_{n;I} \neq \emptyset$. It was conjectured  that  dim($P^{(k,s)}_{n;I})=\mbox{dim}(P^{(k,s)}_{n})$ holds in general provided that $P^{(k,s)}_{n;I} \neq \emptyset$~\cite{Appa2006}.
However, Corollary~\ref{cor:ndiv} suggests that this conjecture may be false for 
$(n,k,s)=(10,6,2)$. (It is not known whether this conjecture is true or false for the $(n,k,s)=(10,6,2)$ case. It is also not known whether 
$P^{(6,2)}_{10;I}\neq \emptyset$.)
Based on the lower bounds for $k$ on  website~\cite{stonybrook}, 
 $P^{(6,2)}_{10;I}$ is the smallest $n,k$ case for $\lambda=1, s=2$ in which this conjecture may fail.
 The following example is consistent with Theorem~\ref{thm:dim} and shows that this conjecture  cannot be generalized as
$\mbox{dim}(P^{(k,s, \lambda)}_{n;I})=\mbox{dim}(P^{(k,s, \lambda)}_{n})$ whenever $P^{(k,s,\lambda)}_{n;I} \neq \emptyset$. 
\begin{ex}
Consider the family of cases $P^{(k,3,\lambda)}_{2;I}$ where 
$8 \lambda/3 \leq k \leq 8\lambda/2$. Theorem 3 in Butler~\cite{Butler2007} implies that for each $\x \in P^{(k,3,\lambda)}_{2;I}$, $J^{\x}_{u}(i_1,\ldots, i_k)=0$ for $u \subseteq [k], \, \forall (i_1,\ldots,i_k)^{\top} \in \{l_1,\ldots,l_n\}^k$ if $|u|$ is odd.
Then $${\rm dim}(P^{(k,3,\lambda)}_{2;I})\leq 2^k-\sum_{j=0}^3{k \choose j}(2-1)^j-\sum_{j=3}^{\lfloor \frac{k+1}{2}\rfloor} {k \choose 2j-1}(2-1)^{2j-1}$$
for $k\in \mathbb{Z}$ such that  $8 \lambda/3 \leq k \leq 8\lambda/2$. On the other hand, for such $k$, assuming that $P^{(k,3,\lambda)}_{2;I}\neq \emptyset$, Theorem~\ref{thm:dim} implies that for odd $\lambda \geq 1$
$${\rm dim}(P^{(k,3,\lambda)}_{2;I})= 2^k-\sum_{j=0}^3{k \choose j}(2-1)^j-\sum_{j=3}^{\lfloor \frac{k+1}{2}\rfloor} {k \choose 2j-1}(2-1)^{2j-1},$$ since
$$ \lambda \binom{3+\ell-1}{\ell-1}\not\equiv 0 \quad {\rm (mod } \ 2)$$ for odd $\ell \in [k-3]$ and odd $\lambda$.
Moreover, Theorem 3 in Butler~\cite{Butler2007}  is consistent with 
Theorem~\ref{thm:dim} as
$$
\lambda \binom{3+\ell-1}{\ell-1}\equiv 0 \quad {\rm (mod } \ 2)$$ for even
 $\ell \in [k-3]$.
\end{ex}
When $n=2$ and  $s$ is even, and for a given feasible point $\x$, the following theorem provides  restrictions for all possible values of $\mbox{dim(Conv}(G(k)^{\rm OD}\,\,\x))$ of ILD~(\ref{ilp:J}) as well as the corresponding sets of equality constraints. These restrictions are stronger than those in Theorem~\ref{thm:dim}.
\begin{thm}\label{thm:dim2}
Let $n=2$, and $s$ be even in ILD~(\ref{ilp:J}). Let $\x$ be a feasible point of ILD~(\ref{ilp:J}),
 $\Omega_1$ and $\Omega'_1$ be as in Theorem~\ref{thm:dimfin}. 
   Then the following hold.
  \begin{enumerate}
  \item[(i)]
   There exists $T\subseteq \Omega_1$ and $T'\subseteq \Omega_1'$ such that
  $${\rm dim(Conv}(G(k)^{\rm OD}\,\, \x))=2^k-\sum_{j=0}^s{k \choose j}-\left(\sum_{d\in T }\left({k \choose d-1}+
  {k \choose d}\right)\right)-\delta_{(|T'|,1)},$$where $k\choose m$ is zero if $m>k$.
  \item[(ii)]
    The equality constraints of 
ILD~(\ref{ilp:J}) together 
with the distinct equalities in 
\begin{align}\label{eqn:convexhullGOD}
	\begin{split}
 J^{\x}_{u''}(i_1,\ldots, i_k)&=0 \quad\forall  u''\subseteq [k] \, \ni\,  |u''| \in T \cup (T-1), 
  \quad \forall \text{$(i_1,\ldots, i_k)^{\top}\in \{l_1,\ldots,l_n\}^k,$}\\
 J^{\x}_{[k]}(i_1,\ldots, i_k)&=0 \quad \text{if $|T'|=|\{k\}|=1$}
 \end{split}
\end{align}
determine  
${\rm Aff}({\rm Conv}(G(k)^{\rm OD}\,\, \x))={\rm Aff}(G(k)^{\rm OD}\,\, \x)$, where the set $T-1$ is obtained from $T$ by subtracting $1$ from each element of $T$.
   \end{enumerate}
\end{thm}
\begin{pf}
The proof follows the proof of Theorem~\ref{thm:dim} up to and including equation~(\ref{eqn:ayni})
 line by line by replacing $S_n \wr S_k$ with $G(k)^{\rm OD}$.
 Now, ${\rm Span}_{\mathbb{Q}}(G(k)^{\rm OD}\,\, \y)$ must be orthogonal to at least one of the irreducible invariant subspaces $W_j$
in the decomposition in Corollary~\ref{cor:GkODQ} for some $j\geq s/2+1$. 
Let $T$ be the set of all even $d$ in $\Omega_1$ 
such that $\mbox{Span}(G(k)^{\rm OD}\,\, \y)$ is orthogonal to $W_{d/2}$.
Let $T'=\{k\}$ if and only if $k$ is odd and 
${\rm Span(} G(k)^{\rm OD}\,\,\y)$ is orthogonal to $W_{(k+1)/2}$.
Then, ($i$) follows since $\mbox{dim}(W_{d/2})={k \choose d-1}+
  {k \choose d}$ if $d \in T$,  $\mbox{dim}(W_{d/2})=\mbox{dim}(W_{k/2})={k \choose k}=1$ if $d=k \in T'$, and  $$\mbox{dim}_{\mathbb{Q}}\mbox{(Span}_{\mathbb{Q}}(G(k)^{\rm OD} \,\,\y))=\mbox{dim(Span}(G(k)^{\rm OD}\,\, \y))=\mbox{dim(Conv}(G(k)^{\rm OD}\,\, \y)).$$   
  By Lemma~\ref{lem:u''2} and equations~(\ref{eqn:mustbezero}),~(\ref{ilp:Jyy}), and~(\ref{eqn:Jxx}) 
\begin{align*}
 J^{\y}_{u''}(i_1,\ldots, i_k)= J^{\x}_{u''}(i_1,\ldots, i_k)=0 &\quad \forall  u''\subseteq[k] \, \ni |u''|\in T\cup (T-1), \,  \forall (i_1,\ldots, i_k)\in \{l_1,\ldots,l_n\}^k,\\
  J^{\y}_{[k]}(i_1,\ldots, i_k)= J^{\x}_{[k]}(i_1,\ldots, i_k)=0 &\quad \text{if } |T'|=|\{k\}|=1, \,  \forall (i_1,\ldots, i_k)\in \{l_1,\ldots,l_n\}^k,
 \end{align*}
   and all feasible points $\x$ of ILD~(\ref{ilp:J}).  
Moreover, the equality constraints of 
ILD~(\ref{ilp:J}) together 
with equations~(\ref{eqn:convexhullGOD}) determine 
$\text{Aff}({\rm Conv}( G(k)^{\rm OD}\,\,\x))=\text{Aff}(G(k)^{\rm OD}\,\,\x)$, proving ($ii$).
\qed
\end{pf}
Next, we prove Theorem~\ref{thm:dimfin}  by  generalizing Theorems~\ref{thm:dim} and~\ref{thm:dim2}. \\ \\
 {\bf Proof of Theorem~\ref{thm:dimfin}}.
 	Let \begin{align*}
 	G=\begin{cases}
 		G(k)^{\rm OD} & \text{if $n=2$ and $s$ is even,}\\
 		G^{\rm iso}(k,n) &\text{otherwise.}
 	\end{cases}
 	\end{align*}
 	Let $\x_1,\ldots,\x_r$ be such that $ 
 	P^{(k,s,\lambda)}_{ n;I}=\text{Conv}\left(\bigcup_{i=1}^rG\,\,\x_i\right)$, and
 	$$\y_i=\x_i-\frac{\lambda n^s}{n^k}\1_{n^k}.$$ 
 	By the proof of Lemma~\ref{lem:dimdim}, $\0_{n^k}\in \text{Conv}\left(\bigcup_{i=1}^r G\,\,\y_i\right)$.
 	Then by Lemma~\ref{lem:dimgengroup}, 
 	{\small
 		\begin{align} \label{eqn:conv2}
 			\text{dim}(P^{(k,s,\lambda)}_{ n;I})=\text{dim(Conv}(\bigcup_{i=1}^rG\,\,\x_i)) =\text{dim(Conv}(\bigcup_{i=1}^rG\,\,\y_i))=\text{dim(Span(}G \,\, \y_1)+\cdots+\text{Span(}G\,\, \y_r)).
 	\end{align}} 
 	Now, we claim that for each $p \in [r]$, 
	\begin{align*} 
		&	\text{dim(Span}(G\,\, \y_1)+\cdots+\text{Span}(G \,\,\y_p))= \\
		&\begin{cases}
			2^k-\sum_{j=0}^s{k \choose j}-\left(\sum_{d\in {T_1}_p }\left({k \choose d-1}+
			{k \choose d}\right)\right)-\delta_{(|{T'_1}_p|,1)}\quad &\text{if $n=2$ and $s$ is even,} \\
			n^k-\sum_{j=0}^s{k \choose j}(n-1)^j-\sum_{d \in {T_2}_p }{k \choose d}(n-1)^{d} & \text{otherwise}
		\end{cases}
	\end{align*} 
	for some ${T'_1}_p\subseteq \Omega'_1$ and ${T_i}_p\subseteq \Omega_i$ such that ${T'_1}_{p}\subseteq {T'_1}_{p-1}\subseteq \Omega'_1$ and ${T_i}_{p}\subseteq {T_i}_{p-1}\subseteq \Omega_i$ for $i=1,2$.
 	By~(\ref{eqn:conv2}), proving the claim, and taking $p=r$ proves Theorem~\ref{thm:dimfin}.
 	
 	First, since  	 
 	\[
 	\text{dim}_{\mathbb{Q}}\text{(Span}_{\mathbb{Q}}(G \,\,\y_1)+\cdots+\text{Span}_{\mathbb{Q}}(G \,\, \y_p))=\text{dim(Span}(G \,\, \y_1)+\cdots+\text{Span}(G \,\, \y_p)), 
 	\] it suffices to prove the claim for $\text{dim}_{\mathbb{Q}}\text{(Span}_{\mathbb{Q}}(G \,\, \y_1)+\cdots+\text{Span}_{\mathbb{Q}}
 	(G \,\,\y_p))$.
 	We prove this claim by induction on $p$. For $p=1$ our claim follows from Theorems~\ref{thm:dim} and~\ref{thm:dim2}.
 	Assume the claim holds for 
 	$p-1$. 
 	Let $U=\text{Span}_{\mathbb{Q}}(G \,\, \y_1)+\cdots+\text{Span}_{\mathbb{Q}}(G \,\, \y_{p-1})$, and 
 	$V=\text{Span}_{\mathbb{Q}}(G \,\, \y_{p})$. Let $\A'\x=\bb'$ be the equality constraints of ILD~(\ref{ilp:J}).
 	Then both $U$ and $V$ are $\mathbb{Q}$-representations of $G$ in Null$_{\mathbb{Q}}(\A')$. By the induction hypothesis the claim holds for 
 	$U$ with ${T'_1}_{p-1}$, ${T_i}_{p-1}$ for $i=1,2$. 
 	Let $T_1, T_1'$ be $T,T'$, in the proof of Theorem~\ref{thm:dim2}, and $T_2$ be $T$  in the proof of Theorem~\ref{thm:dim} for $V=\text{Span}_{\mathbb{Q}}(G \,\, \y_{p})$ in their corresponding cases. Let ${T_i}_p=T_{ip-1}\cap T_i$ for $i=1,2$, and ${T_1}_p'=T_{1p-1}'\cap T_1'$. Then, ${T'_1}_{p}\subseteq {T'_1}_{p-1}\subseteq \Omega'_1$, ${T_i}_{p}\subseteq {T_i}_{p-1}\subseteq \Omega_i$ for $i=1,2$, and
 	\begin{align*} 
 		\text{dim}_{\mathbb{Q}}(U+V)=\begin{cases}
 			2^k-\sum_{j=0}^s{k \choose j}-\left(\sum_{d\in {T_1}_p }\left({k \choose d-1}+
 			{k \choose d}\right)\right)-\delta_{(|{T'_1}_p|,1)} &\text{if $n=2$ and $s$ is even,}  \\
 			n^k-\sum_{j=0}^s{k \choose j}(n-1)^j-\sum_{d\in {T_2}_p }{k \choose d}(n-1)^{d} &\text{otherwise.} \quad \quad \quad \quad \quad \, \qed
 		\end{cases} 
 	\end{align*} 
 
 The following corollary follows from the proof of Theorem~\ref{thm:dimfin}.
\begin{cor}
 The equality constraints of  ILD~(\ref{ilp:J}) together 
with equations~(\ref{eqn:convexhullGOD}) if $n=2$ and $s$ is even, and with equations~(\ref{eqn:convexhullGiso}) if either $n\geq 3$ or $s$ is odd,  determine the affine hull of the 
feasible points of ILD~(\ref{ilp:J}) 
or equivalently ILD~(\ref{ilp:BF}).
\end{cor}
\section{Generalization to ILDs with equality constraints and discussion}\label{sec:generalization}
A feasible LD with no redundant  constraints  
and  no inequalities  satisfied by all feasible $\x$ as an equality is said to be in {\em standard form}.
Since the feasible set of any feasible LD can be made the feasible set of an LD in standard form~\cite{Geyer2018}, WLOG let LD~(\ref{eqn:geneqILP}) be in  {\em standard form}. 
Let $P$  be the feasible set of LD~(\ref{eqn:geneqILP}),   
and $P_I$ be the convex hull of the feasible points of ILD~(\ref{eqn:geneqILP}). Method 4 in~\cite{Geyer2018} can be used for finding $G^{\text{LD}(\ref{eqn:geneqILP})}$.

In this section, we discuss how in general  many  possible values for the dimension of the convex hull of all feasible points of an ILD with the LD relaxation symmetry group 
$G^{\text{LD}(\ref{eqn:geneqILP})}$ can be ruled out. We also describe how the zero right hand side linear equality constraints associated with $G^{\text{LD}(\ref{eqn:geneqILP})}$ can be generated. These  are the only sets of zero right hand side linear equality constraints associated with $G^{\text{LD}(\ref{eqn:geneqILP})}$ up to equivalence that can potentially be implied by the integrality constraints of the ILD.
 All the  results of this section are valid if $G^{\text{LD}(\ref{eqn:geneqILP})}$ is replaced with any other subgroup $G$ of the symmetry group of 
ILD~(\ref{eqn:geneqILP}) provided that ${\rm Row}(\A)$ and
 ${\rm Row}(\A)^{\perp}$  are both $G$-invariant subspaces (i.e.,
 subrepresentations). This requirement necessarily holds for
$G = G^{{\rm LD}(\ref{eqn:geneqILP})}$ if LD~(\ref{eqn:geneqILP}) is in standard
form. We use $G^{\text{LD}(\ref{eqn:geneqILP})}$, as it is the largest known subgroup of the symmetry group of 
ILD~(\ref{eqn:geneqILP}) for which there is a known 
 generation method~\cite{Geyer2018} without finding all solutions.

Let $\mathbb{F}^Y$ be the vector space of vectors indexed by the index set $Y$ of variables of ILD~(\ref{eqn:geneqILP}).
Let $G$ be a subgroup of the group of all permutations of the elements of $Y$. Let 
\begin{equation*}\label{eqn:RgF}
R:G\rightarrow \rm{Aut}_{\mathbb{F}}(\mathbb{F}^{Y})
\end{equation*}
be the permutation $\mathbb{F}$-representation associated with $Y$, where $G$ acts on $\mathbb{F}^Y$ by $R(g)f(y)=f(g^{-1}y)$ for all $f \in \mathbb{F}^Y$.
Then by Maschke's theorem (cf.~\cite{Goodman}, Theorem 2.4.1), 
\begin{align}\label{eqn:decomp}	
	R=\bigobot_{i=1}^{b_{\mathbb{F}}} R_i:G&\rightarrow \text{GL}\left(V_1^{\mathbb{F}} \obot \cdots \obot V_{b_{\mathbb{F}} }^{\mathbb{F}}\right), \nonumber\\
\mathbb{F}^Y&=V_1^{\mathbb{F}} \obot \cdots \obot V_{b_{\mathbb{F}} }^{\mathbb{F}},
\end{align}
where each $(R_i,V_i^{\mathbb{F}})$ is an irreducible $\mathbb{F}$-representation of $G$. For a coordinate vector $\x \in \F^Y$, we denote $R(g)\x$ by $g\x$, where we view $g \in G$ as a permutation or a permutation matrix depending on the context. The direct sum in  decomposition~(\ref{eqn:decomp})  can be taken to be an orthogonal direct sum with respect to the complex dot product as $G$ acts on $\mathbb{F}^Y$ with unitary matrices, and unitary matrices preserve the complex dot product $\langle \cdot \mid \cdot \rangle$, i.e.,  
$$\langle g\x\mid g\y\rangle=
\langle \x\mid \y\rangle\ \forall g \in G, \text{ and } \x,\y \in \mathbb{F}^Y.$$  
 %

For a representation $(\rho,V)$ of $G$ let 
$${\rm Fix}_{V} (G)=\{\vv \in V\mid \rho(g)\vv =\vv\} $$
be the fixed subspace of $V$ under the action of $G$.
Then, we have the following theorem.
\begin{thm}\label{thm:Reviewer3}
Let $S$ be the feasible set of ILD~(\ref{eqn:geneqILP}), and $\F = \Q$, $\R$, or $\C$. Then ${\rm Aff}_\F (S) = \p + U$, where 
$\p \in {\rm Conv}(S) \cap {\rm Fix}_{\Q^n} (G)$, and 
$U \subseteq {\rm Row}_{\F} (\A)^{\perp}$ is a subrepresentation for $G=G^{{\rm LD}(1)}$.
\end{thm}
\begin{pf}
	First, $\E=(1/|G|) \sum_{g \in G}
	g$ is the orthogonal projection matrix onto the fixed space ${\rm Fix}_{\F^n} (G)$. 
 For $\x_0 \in S$, $G\x_0 \subseteq S$. Then	
$$ \p :=
\frac{1}{|G|} \sum_{g \in G}
g\x_0 \in {\rm Conv}(S) \cap {\rm Fix}_{\Q^n}(G) \subset {\rm Aff}_\F (S) \cap {\rm Fix}_{\Q^n} (G).$$
  Let
$U ={\rm Span}_\F \{\sss - \p \mid \sss \in S\}$. Then $U$ is a subrepresentation, since $\{\sss - \p \mid \sss \in S\}$ is a $G$-invariant
subset. Then, ${\rm Aff}_{\F} (S) = \p + U$ and
 $U \subseteq {\rm Row}_\F (\A)^{\perp}$, because $\A\sss = \bbb$ for feasible $\sss$ and  $\A\p = \bbb$.
 \qed
\end{pf}
\begin{cor}\label{cor:referee}
There is a $G^{{\rm LD} (\ref{eqn:geneqILP})}$ invariant subspace $U \subseteq {\rm Row}_\Q (\A)^{\perp}$ such that ${\rm dim}( P_I )= {\rm dim}_\Q(U ).$
\end{cor}
The next theorem is a generalization of Theorem~\ref{thm:dim}~$(i)$ and follows from Corollary~\ref{cor:referee}.
\begin{thm}\label{thm:nointeger}
 Let
\begin{equation}\label{eqn:refree}
	{\rm Row}_\Q (\A)^\perp = m_1 W_1 \obot\cdots\obot m_{\ell} W_{\ell}
\end{equation}
be a decomposition into irreducible subrepresentations of $G^{{\rm LD} (\ref{eqn:geneqILP})}$, where the different $W_j$s are non-equivalent and $W_j$ occurs with multiplicity $m_j$ with $m_j\geq 1$. Then $${\rm dim}(P_I)\in \left\{ \sum_{j=1}^{\ell}
m'_j {\rm dim}_\Q (W_j )\mid 0 \leq m'_j \leq m_j \text{ for } j \in [\ell]
\right\}.$$
\end{thm}
\begin{pf}
 Every $\Q$-subrepresentation $U$ of ${\rm Row}_\Q (\A)^{\perp}$ is equivalent
to one of the form
$$m'_1 W_1 \obot \cdots \obot m'_\ell W _\ell$$
 such that   $0 \leq m'_j \leq m_j$ for $j \in [\ell]$, and thus by Corollary~\ref{cor:referee} 
\begin{align}
{\rm dim}(P_I)= {\rm dim}_\Q(U )=\left\{ \sum_{j=1}^{\ell}
m'_j {\rm dim}_\Q (W_j )\mid 0 \leq m'_j \leq m_j
\text{ for } j \in [\ell]\right\}. \qqed
\end{align}
\end{pf}
For the orthogonal array polytope $P^{(k,s,\lambda)}_{n;I}$, $\p=\lambda n^{s-k}\1_{n^k}$ in Lemma~\ref{lem:dimdim},
$G(k)^{{\rm OD}} \leq G^{{\rm LD}(1)}$ if $n=2$ and $s$ is even, $S_n \wr S_k\leq G^{{\rm LD}(1)}$ otherwise, and $m_j = 1$ for $j \in [\ell]$. Moreover,  
$(W_j )_\C$, i.e., the $\CCC$-representation obtained from $W_j$ by extending the field of scalars of $W_j$  to $\CCC$ is still irreducible.
In general, $(W_j )_\C$ may not be irreducible~\cite[Theorem 9.21]{Isaacs1994}. Fortunately, there is a randomized algorithm 
in~\cite{Babai} that runs in expected polynomial time for computing a decomposition of 	${\rm Row}_\C (\A)^\perp$ into irreducible subrepresentations.
On the other hand, developing an algorithm that finds a decomposition as in decomposition~(\ref{eqn:refree}) in polynomial time is a problem in representation theory proposed 
as an open problem  by Babai and R\'{o}nyai~\cite[Problem 7.1]{Babai}.
 Unfortunately, knowing a decomposition of ${\rm Row}_\C (\A)^\perp$ into irreducible $\C$-subrepresentations yields 
a weaker result compared to  Theorem~\ref{thm:nointeger}.

Let
 \begin{equation}\label{eqn:complexification}
 {\rm Row}_\C (\A)^\perp=(\text{Row}_{\mathbb Q} (\A)^\perp)_{\C}=\bigobot_{i}c_iV_i
 \end{equation} 
be a  decomposition 
into irreducible 
$\C$-subrepresentations, where $c_i \geq 1$. Then the projection matrix onto each subspace $c_iV_i$ is known~\cite[Theorem 8]{Serre}. Hence, by the uniqueness of decomposition~(\ref{eqn:complexification}) into subspaces $c_iV_i$~\cite[Theorem 8]{Serre}, these projection matrices are  necessarily orthogonal as the underlying representation is unitary. Then, the orthogonal projection matrix onto each $m_jW_j$ in decomposition~(\ref{eqn:refree}) is also known as 
each $m_jW_j$ is an orthogonal direct sum of some $c_iV_i$s~\cite[Theorem 9.21]{Isaacs1994}. Hence,  if $c_i=1$ for each $i$, then $m_j=1$ for each $j \in [\ell]$, and 
 it is easy to obtain the orthogonal projection matrix onto $W_j$ for $j \in [\ell]$.

Let $\mathbb{A}$ and $\mathbb{A}^{{\rm LD}(\ref{eqn:geneqILP})}$
be the affine spaces where the convex hull of all feasible points of  ILD~(\ref{eqn:geneqILP}) and 
 LD~(\ref{eqn:geneqILP}) lie.
Then dim$(\mathbb{A})$  may be smaller than 
dim$(\mathbb{A}^{{\rm LD}(\ref{eqn:geneqILP})})$ due to the integrality constraints.  It is far from clear what additional equality 
constraints are needed to obtain  $\mathbb{A}$.
For cases in which a large group of permutations preserves the feasible set of the ILD, the representation theory based approach in this paper provides a method to obtain a  finite collection of candidate sets of equality constraints  that correspond to a finite  set of candidate affine subspaces for  $\mathbb{A}$. 
In particular, if dim$(\mathbb{A})<\mbox{dim}(\mathbb{A}^{{\rm LD}(\ref{eqn:geneqILP})})$,  let $$W=\{\vv \in \text{Row}_{\mathbb Q} (\A)^\perp \mid \vv^{\top}\x = cons(\vv) \, \,  \forall  \x \text{ that is a feasible point of ILD~(\ref{eqn:geneqILP})}\},$$ where $cons(\vv) \in \mathbb{Q}$ is a constant that depends on $\vv$. Then $W$ is a subspace of $\text{Row}_{\mathbb Q} (\A)^\perp$. Moreover, for each $\vv \in W$, 
$$ \vv^{\top}\x=cons(\vv) \implies \vv^{\top}g\x=(g^{\top}\vv)^{\top}\x=cons(\vv) \, \, \forall g \in G^{{\rm LD}(\ref{eqn:geneqILP})}.$$ Hence, $W$ is a 
$\mathbb{Q}$-representation of $G^{{\rm LD}(\ref{eqn:geneqILP})}$, and consequently there exists a collection of mutually non-equivalent irreducible $\mathbb{Q}$-representations 
$W_{i_1}, \ldots, W_{i_{l'}} $ of $G^{{\rm LD}(\ref{eqn:geneqILP})}$ in $\text{Row}_{\mathbb Q} (\A)^\perp$ such that 
\begin{equation}\label{eqn:Wij}
W=m''_{i_1}W_{i_1}\obot\cdots\obot m''_{i_{l'}}W_{i_{l'}}.
\end{equation}
 
 Let $d_{i,k} \in \mathbb{Q}$ be such that
$(g^{\top}\vv_{i,k})^{\top}\x=d_{i,k}$  $\forall \x \in \mathbb{A}$ and $g \in G^{{\rm LD}(\ref{eqn:geneqILP})}$, $i \in [r_k]$, $k \in [\sum_jm''_{i_j}]$, 
 where $\{\vv_{1,k},\ldots,\vv_{r_k,k}\}$ is a basis for the $k$th irreducible subrepresentation in  decomposition~(\ref{eqn:Wij}). This implies 
\begin{equation}\label{eqn:dik}
\left (\vv_{i,k}-\frac{1}{|G^{{\rm LD}(\ref{eqn:geneqILP})}|}\sum_{g \in G^{{\rm LD}(\ref{eqn:geneqILP})}}g^{\top}\vv_{i,k}\right)^{\top}\x=0
\end{equation}
 $\forall \x \in \mathbb{A}$ and  $i\in [r_k]$, $k \in [\sum_jm''_{i_j}]$. Hence, by using representation theory it is possible to generate candidate constraints satisfied by every point of $\mathbb{A}$
as the zero right hand side linear equality constraints associated with $G^{\text{LD}(\ref{eqn:geneqILP})}$. 
In particular, if $c_i=1$ for all $i$ in decomposition~(\ref{eqn:complexification}), it is easy to generate 
such constraints as by the above discussion, it is easy to find the bases $\{\vv_{1,k},\ldots,\vv_{r_k,k}\}$ for $k\in [\ell']$. 
When the goal is to find a solution instead of finding all non-isomorphic solutions with respect to $G^{{\rm LD}(\ref{eqn:geneqILP})}$, $\mathbb{A}$ can be assumed  to be the affine space where the convex hull of the orbit of one solution $\x$ under the action of $G^{{\rm LD}(\ref{eqn:geneqILP})}$ (the isomorphism class of
 $\x$ with respect to $G^{{\rm LD}(\ref{eqn:geneqILP})}$) lie. For such an 
 $\mathbb{A}$, potentially  dim$(\mathbb{A})<\mbox{dim}(\mathbb{A}^{{\rm LD}(\ref{eqn:geneqILP})})$ making it possible to find solutions
  after incorporating the constraints~(\ref{eqn:dik})  for each irreducible $\mathbb{Q}$-representation in some collection of irreducible $\mathbb{Q}$-representations $W_{i_1}, \ldots, W_{i_{l'}}$ of $G^{{\rm LD}(\ref{eqn:geneqILP})}$ 
all in $\text{Row}_{\mathbb Q} (\A)^\perp$. The hope is that the additional constraints would render the resulting ILD to be easier to solve, where proving infeasibility or finding a solution can be accomplished by using the altered version of the isomorphism pruning algorithm of~\cite{Margot2007} as in~\cite{Bulutoglu2008,Bulutoglu2016} after converting the ILD to an ILP by introducing the zero objective function. 
 Finally, one can iterate over many different collections of irreducible subrepresentations, and try solving several ILDs until a  solution is found.  If a larger subgroup $H$ of the symmetry group of the ILD containing  $G^{{\rm LD}(\ref{eqn:geneqILP})}$ is used, then the number of irreducible subrepresentations in decomposition~(\ref{eqn:refree}) will be decreased. This will not only decrease the number of ILDs that need to be solved, but also potentially decrease the difficulty of the resulting ILPs obtained from the resulting ILDs due to having additional constraints.
Hence, this method will be most useful for finding a feasible point to  an ILD 
 for which a large subgroup of its symmetry group is known and
finding a feasible point is  computationally challenging.

 In this article we reveal the underlying representation theory that 
dictates the results regarding dim($P^{(k,s)}_{n;I})$ 
in~\cite{Appa2006,AppaJanssen2006,Balas1989,Balinski1974,Euler1987,
Euler1986}.
For $P^{(k,s,\lambda)}_{n;I}\neq \emptyset$, we not only provide a sufficient condition  for   
\begin{equation}\label{eqn:finconj}
\mbox{dim}(P^{(k,s,\lambda)}_{n;I})=\mbox{dim}(P^{(k,s,\lambda)}_{n})=n^k-\sum_{j=0}^s{k \choose j}(n-1)^j
\end{equation}
to be true, we also provide a family of examples with $P^{(k,s,\lambda)}_{n;I}\neq \emptyset$ such that equation~(\ref{eqn:finconj}) is not valid  when this sufficient condition is not  satisfied. 
We develop our method of proof into  Theorem~\ref{thm:nointeger} that finds restrictions for the dimension of the affine hull of the feasible set of an arbitrary ILD~(\ref{eqn:geneqILP}) with LD relaxation symmetry group  $G^{\text{LD}(\ref{eqn:geneqILP})}$. Additionally, we discuss how to determine the sets of potential zero right hand side  linear equality constraints associated with $G^{\text{LD}(\ref{eqn:geneqILP})}$ that come together with the restrictions from Theorem~\ref{thm:nointeger}.
 Based on Theorem~\ref{thm:nointeger}, we then propose a heuristic for finding a feasible point to an ILD for which a large subgroup of its symmetry group is known and
finding a feasible point is  computationally challenging. 
Theorem~\ref{thm:nointeger} requires knowing  decomposition~(\ref{eqn:decomp}) for $\mathbb{F}=\mathbb{Q}$. Developing an algorithm for determining  decomposition~(\ref{eqn:decomp}) for $\mathbb{F}=\mathbb{Q}$ in polynomial time is a problem in representation theory proposed 
as an open problem (Problem 7.1) by Babai and R\'{o}nyai~\cite{Babai}.
We emphasize the applicability of a solution to this problem to determining the feasibility of an ILP for which there is a large known subgroup of symmetries. 
\section*{Acknowledgments}
The author thanks Dr. William P. Baker for solving a partial difference equation for deriving equation~(\ref{eqn:formula}) and 2nd Lt Kristopher Kilpatrick for a careful reading of the paper that lead to several improvements. The author also thanks an Associate Editor and three anonymous referees for their helpful comments that greatly improved the paper.

  The views expressed in this article are those of the author, and do not reflect the official policy or position of the United States Air Force, Department of Defense, or the U.S.~Government.
  
\bibliographystyle{elsarticle-harv}
\bibliography{BibliographyDO}

\end{document}